\documentclass[reqno]{amsart}
\usepackage{amsfonts,amsmath,amssymb}

\numberwithin{equation}{section}

\newtheorem{thm}{Theorem}[section]
\newtheorem{lem}{Lemma}[section]

\newtheorem{Step}{Step}[section]

\begin{document}
\title[Sharp Sobolev inequalities]
{Sharp Sobolev inequalities for\\
vector valued maps}
\author{Emmanuel Hebey}
\address{Emmanuel Hebey, Universit\'e de Cergy-Pontoise, 
D\'epartement de Math\'ematiques, Site de 
Saint-Martin, 2 avenue Adolphe Chauvin, 
95302 Cergy-Pontoise cedex, 
France}
\email{Emmanuel.Hebey@math.u-cergy.fr}

\date{July 18, 2005}

\maketitle

Let $(M,g)$ be a smooth compact Riemannian $n$-manifold, $n \ge 3$. Let also $p \ge 1$ be an integer, 
and $M_p^s({\mathbb R})$ be the vector space of symmetrical $p\times p$ real matrix. 
Namely the vector space of $p\times p$ real matrix $S = (S_{ij})$ which are such that $S_{ij} = S_{ji}$ 
for all $i, j = 1,\dots,p$. 
We regard the elements $S = (S_{ij})$ in $M_p^s({\mathbb R})$ as bilinear forms on ${\mathbb R}^p$ by letting 
$S(X,Y) = \sum_{i,j=1}^pS_{ij}X^iY^j$, where $X = (X^1,\dots,X^p)$ and $Y = (Y^1,\dots,Y^p)$. 
We let $H_{1,p}^2(M)$ be the Sobolev space consisting of $p$-maps 
${\mathcal U}: M \to {\mathbb R}^p$, ${\mathcal U} = (u_1,\dots,u_p)$, which are 
such that the $u_i$'s are all in the standard Sobolev space $H_1^2(M)$ consisting of functions in $L^2(M)$ with 
one derivative in $L^2$. We let also 
$A: M \to M_p^s({\mathbb R})$ 
smooth, $A = (A_{ij})$, be such that $A(x)$ is positive for all $x \in M$ as a bilinear form. 
Then we consider Sobolev inequalities like
\begin{equation}\label{GenSobIneqIntro}
\left(\int_M\vert{\mathcal U}\vert^{2^\star}dv_g\right)^{2/2^\star} 
\le K \int_M\vert\nabla{\mathcal U}\vert^2dv_g + \Lambda\int_MA({\mathcal U},{\mathcal U})dv_g
\hskip.1cm ,
\end{equation}
where the inequality is required to hold for all ${\mathcal U} \in H_{1,p}^2(M)$, $K, \Lambda > 0$ are positive constants, $A$ is regarded as a bilinear form, 
$dv_g$ is the Riemannian volume element of $g$, 
and the exponent $2^\star = 2n/(n-2)$ is critical from the Sobolev viewpoint. 
Here, in (\ref{GenSobIneqIntro}), 
$\vert{\mathcal U}\vert^{2^\star} = \sum_{i=1}^p\vert u_i\vert^{2^\star}$, 
$\vert\nabla{\mathcal U}\vert^2 = \sum_{i=1}^p\vert\nabla u_i\vert^2$, and 
$A({\mathcal U},{\mathcal U}) = \sum_{i,j=1}^pA_{ij}u_iu_j$ when we write that ${\mathcal U} = (u_1,\dots,u_p)$. 
The questions we ask in this paper are: what is the value $K_s$ of the sharp $K$ in (\ref{GenSobIneqIntro}), 
does the corresponding sharp inequality hold true, and if yes, 
does its saturated version (where $\Lambda$ is lowered to its minimum 
value under the constraint $K = K_s$) possess extremal functions. When 
$p = 1$, we are back to the classical setting of the Sobolev inequality for functions. 
Inequality (\ref{GenSobIneqIntro}) when $p \ge 2$, and the above sequence of questions, are the natural 
extensions to vector valued maps of the $AB$-program which was developed 
in the case of functions. 
Possible references in 
book form for the problem in the case of functions and the $AB$-program are Druet and Hebey \cite{DruHebAMSBook}, and Hebey 
\cite{HebCIMSBook}. In what follows we let $K_n$ be the sharp constant $K$ in the Euclidean Sobolev inequality
$\Vert u\Vert_{2^\star} \le K_n\Vert\nabla u\Vert_2$. Then
\begin{equation}\label{SharpCstEucl}
K_n = \sqrt{\frac{4}{n(n-2)\omega_n^{2/n}}}\hskip.1cm ,
\end{equation}
where $\omega_n$ is the volume of the unit $n$-dimensional sphere. In the sequel, 
we say that a matrix $A = (A_{ij})$ is 
{\it cooperative} if $A_{ij} \ge 0$ for all $i \not= j$. When $A: M \to M_p^s({\mathbb R})$ is a map, 
$A$ is said to be 
cooperative in $M$ if $A_{ij}(x) \ge 0$ for all $i \not= j$, and all $x \in M$. We also say that $A$ 
is {\it globally irreducible} if the index set $\{1,\dots,p\}$ does not split in two disjoint subsets 
$\{i_1,\dots,i_k\}$ and $\{j_1,\dots,j_{k^\prime}\}$, $k + k^\prime = p$, such that $A_{i_\alpha j_\beta} \equiv 0$ 
for all $\alpha = 1,\dots,k$ and $\beta = 1,\dots,k^\prime$. 
A $p$-map 
${\mathcal U} = (u_1,\dots,u_p)$ is said to be 
{\it nonnegative} if the $u_i$'s are all nonnegative functions (i.e $u_i \ge 0$ for all $i$), 
{\it weakly positive} if the $u_i$'s are all positive functions unless they are identically zero 
(i.e, for any $i$, either $u_i > 0$ or $u_i \equiv 0$), and 
{\it strongly positive} if the $u_i$'s are all positive functions (i.e $u_i > 0$ for all $i$). The map 
${\mathcal U}$ is said to be 
of {\it undeterminate sign} if neither ${\mathcal U}$ nor $-{\mathcal U}$ is nonnegative. 
For $0 < \theta < 1$, we let $C^{2,\theta}_p(M)$ be the space of $p$-maps with components 
in $C^{2,\theta}(M)$. If ${\mathcal H}$ is a subset of $C^{2,\theta}_p(M)$, invariant under the scaling 
${\mathcal U} \to \lambda{\mathcal U}$ for $\lambda$ real, 
we refer to the $L^{2^\star}$-{\it normalized subset} of ${\mathcal H}$ as the subset of 
${\mathcal H}$ consisting of the ${\mathcal U}$ in ${\mathcal H}$ such that 
$\int_M\vert{\mathcal U}\vert^{2^\star}dv_g = 1$.
At last, we let $S_g$ be the scalar curvature of $g$. Our result states as follows.

\begin{thm}\label{TheThm} Let $(M,g)$ be a smooth compact Riemannian manifold 
of dimension $n \ge 3$, $p \ge 1$ an integer, 
and $A: M \to M_p^s({\mathbb R})$ 
smooth and such that $A(x)$ is positive as a bilinear form for all $x \in M$.
The value $K_s$ 
of the sharp constant $K$ in (\ref{GenSobIneqIntro}) is
$K_s = K_n^2$, where $K_n$ is given by (\ref{SharpCstEucl}), 
and there exists $\Lambda > 0$ such that the sharp inequality
\begin{equation}\label{SharpIneq}
\left(\int_M\vert{\mathcal U}\vert^{2^\star}dv_g\right)^{2/2^\star} 
\le K_n^2 \int_M\vert\nabla{\mathcal U}\vert^2dv_g +\Lambda \int_MA({\mathcal U},{\mathcal U})dv_g
\end{equation}
holds for all ${\mathcal U} \in H_{1,p}^2(M)$, 
where $A$ is regarded as a bilinear form. Moreover, if 
$\Lambda_0(g)$ stands for the lowest $\Lambda$ in (\ref{SharpIneq}), then $\Lambda_0(g) > 0$ and, 
when $n \ge 4$,
\begin{equation}\label{2ndCstIneq}
A_{ii}(x) \ge \frac{(n-2)K_n^2}{4(n-1)\Lambda_0(g)} S_g(x)
\end{equation}
for all $i = 1,\dots,p$ and all $x \in M$, where the $A_{ij}$'s are the components of 
$A$. At last, 
if the inequality in (\ref{2ndCstIneq}) is strict for all $i$ and all $x$, then the sharp and saturated inequality
\begin{equation}\label{SharpSatIneq}
\left(\int_M\vert{\mathcal U}\vert^{2^\star}dv_g\right)^{2/2^\star} 
\le K_n^2 \int_M\vert\nabla{\mathcal U}\vert^2dv_g + 
\Lambda_0(g)\int_MA({\mathcal U},{\mathcal U})dv_g
\end{equation}
possesses extremal maps, namely nontrivial $p$-maps in $C^{2,\theta}_p(M)$ which realize the 
equality in (\ref{SharpSatIneq}), and the 
$L^{2^\star}$-normalized set of such extremal maps is precompact in the 
$C^{2,\theta}_p$-topology, where $0 < \theta < 1$. When no specific assumption is made on $A$, extremal 
maps might be of undeterminate sign. They can be 
chosen weakly positive if $-A$ is cooperative, and strongly positive if $A$ is also globally irreducible.
\end{thm}

When $p = 1$, as already mentionned, 
we are back to the classical Sobolev inequality for functions. The validity of the classical 
sharp inequality for functions 
on arbitrary manifolds was proved in Hebey and Vaugon  \cite{HebVau1,HebVau2}. The existence 
of extremal functions (and the above result when $p=1$) was
studied in Djadli and Druet \cite{DjaDru}, and Hebey and Vaugon \cite{HebVauMathZ}. Possible references in
book form on the sharp classical Sobolev inequality 
are Druet and Hebey \cite{DruHebAMSBook}, and Hebey \cite{HebCIMSBook}. We refer also to 
Collion, Hebey and Vaugon \cite{ColHebVau}, Ghoussoub and Robert \cite{GhoRob}, Humbert and 
Vaugon \cite{HumVau}, and Robert \cite{Rob}. An easy corollary to Theorem \ref{TheThm} is that 
sharp and saturated inequalities like (\ref{SharpSatIneq}) 
always possess extremal maps when $(M,g)$ has nonpositive scalar curvature and $n \ge 4$. Another possible corollary, 
which follows from the theorem and the resolution of the Yamabe problem by 
Aubin \cite{Aub} and Schoen \cite{Sch}, see Section \ref{SectProof3}, is that 
if $n \ge 4$, $A$ does not depend on $x$, and $(M,g)$ has constant scalar curvature, then (\ref{SharpSatIneq}) 
possesses extremal maps. 

\medskip We prove Theorem \ref{TheThm} in Sections 
\ref{SectProof1}, \ref{SectProof2}, and \ref{SectProof4} below. That $K_s = K_n^2$, and that 
(\ref{SharpIneq}) and (\ref{2ndCstIneq}) are true, easily follows from what has been done in the scalar case 
of Sobolev type inequalities. This is discussed in Section \ref{SectProof1}.
The difficult part is to prove 
existence and compactness of extremal maps. 
This is discussed in Sections \ref{SectProof2} and \ref{SectProof4}. Section  \ref{SectProof3} 
is devoted to the proof of the above mentionned corollaries to Theorem \ref{TheThm}. 
When $p \ge 2$, 
contrary to the scalar case where the maximum principle 
for functions can be applied, there are no maximum principle for the  
equations associated to 
inequalities like (\ref{SharpIneq}). We have to deal with maps of undeterminate sign, and not 
only with positive, or even nonnegative maps. The case of maps is more involved than the case of 
functions. 

\section{Proof of the first part of Theorem \ref{TheThm}}\label{SectProof1}

We let $(M,g)$ be a smooth compact Riemannian manifold 
of dimension $n \ge 3$, $p \ge 1$ an integer, 
and $A: M \to M_p^s({\mathbb R})$, $A = (A_{ij})$, be 
smooth and such that $A(x)$ is positive as a bilinear form for all $x \in M$. 
We know from Hebey and Vaugon \cite{HebVau1,HebVau2} that there exists $B > 0$ such that for any 
$u \in H_1^2(M)$,
\begin{equation}\label{ProofEqt1}
\left(\int_M\vert u\vert^{2^\star}dv_g\right)^{2/2^\star} \le 
K_n^2\int_M\vert\nabla u\vert^2dv_g + B\int_Mu^2dv_g\hskip.1cm ,
\end{equation}
where $H_1^2(M)$ is the Sobolev space of functions in $L^2(M)$ with one derivative in $L^2$. 
Since $2/2^\star \le 1$, $(a + b)^{2/2^\star} \le a^{2/2^\star} + b^{2/2^\star}$ for $a, b \ge 0$, and it 
follows from (\ref{ProofEqt1}) that for any ${\mathcal U} \in H_{1,p}^2(M)$,
\begin{equation}\label{ProofEqt2}
\left(\int_M\vert{\mathcal U}\vert^{2^\star}dv_g\right)^{2/2^\star} 
\le K_n^2 \int_M\vert\nabla{\mathcal U}\vert^2dv_g + B\int_M\vert{\mathcal U}\vert^2dv_g
\hskip.1cm ,
\end{equation}
where $H_{1,p}^2(M)$ is the space of $p$-maps 
${\mathcal U}: M \to {\mathbb R}^p$, ${\mathcal U} = (u_1,\dots,u_p)$, which are 
such that the $u_i$'s are all in $H_1^2(M)$. Since we assumed that $A(x)$ is positive for all $x$ as 
a bilinear form, there exists $t > 0$ such that $\delta_{ij} \le tA_{ij}(x)$ for all $x$, in the sense of 
bilinear forms. Letting $\Lambda = Bt$, we get that for any ${\mathcal U} \in H_{1,p}^2(M)$,
\begin{equation}\label{ProofEqt3}
\left(\int_M\vert{\mathcal U}\vert^{2^\star}dv_g\right)^{2/2^\star} 
\le K_n^2 \int_M\vert\nabla{\mathcal U}\vert^2dv_g + \int_MA_\Lambda({\mathcal U},{\mathcal U})dv_g
\hskip.1cm ,
\end{equation}
where $A_\Lambda = \Lambda A$, and $A$ is regarded as a bilinear form. Conversely, 
taking ${\mathcal U}$ in (\ref{GenSobIneqIntro}) such that the components $u_j$ 
of ${\mathcal U}$ are all zero if $j \not= i$, and $u_i = u$ is arbitrary, 
it clearly follows from 
an inequality like (\ref{GenSobIneqIntro}) that for any $i$, and any $u \in H_1^2(M)$,
\begin{equation}\label{ProofEqt4}
\left(\int_M\vert u\vert^{2^\star}dv_g\right)^{2/2^\star} \le 
K\int_M\vert\nabla u\vert^2dv_g + \Lambda\int_MA_{ii}u^2dv_g\hskip.1cm .
\end{equation}
In particular, see for instance Hebey \cite{HebCIMSBook}, we get from 
(\ref{ProofEqt4}) that we necessarily have that 
$K \ge K_n^2$. This inequality, together with (\ref{ProofEqt3}), gives  
that $K_s = K_n^2$. By (\ref{ProofEqt3}) we also have that (\ref{SharpIneq}) is true, and 
by the definition of $\Lambda_0(g)$, we can take 
$\Lambda = \Lambda_0(g)$ in (\ref{SharpIneq}). In particular, 
(\ref{SharpSatIneq}) is true. Like when passing 
from (\ref{GenSobIneqIntro}) to (\ref{ProofEqt4}), 
it follows from the sharp and saturated (\ref{SharpSatIneq}) that 
for any $i$, and any $u \in H_1^2(M)$,
\begin{equation}\label{ProofEqt5}
\left(\int_M\vert u\vert^{2^\star}dv_g\right)^{2/2^\star} \le 
K_n^2\int_M\vert\nabla u\vert^2dv_g + \Lambda_0(g)\int_MA_{ii}u^2dv_g\hskip.1cm .
\end{equation}
Taking $u = 1$ in (\ref{ProofEqt5}), we get that $\Lambda_0(g)\int_MA_{ii}dv_g \ge V_g^{2/2^\star}$ 
for all $i$, where $V_g$ is the volume of $M$ with respect to $g$. Since $A(x)$ is positive for all $x$, 
the $A_{ii}$'s are positive functions, and $\Lambda_0(g)$ has to be positive.
By the developments in Aubin \cite{Aub}, 
we also get with (\ref{ProofEqt5}) that when $n \ge 4$, (\ref{2ndCstIneq}) has to be true 
for all $x \in M$, and all $i$. More precisely, given $\delta > 0$ small, 
$\varepsilon > 0$ small, and $x_0 \in M$, we let $u_{x_0}^\varepsilon$ be the function defined by
$$u_{x_0}^\varepsilon = \left(\varepsilon + r^2\right)^{1-n/2} 
- \left(\varepsilon + \delta^2\right)^{1-n/2}$$
if $r \le \delta$, and $u_{x_0}^\varepsilon = 0$ if not, where $r = d_g(x_0,\cdot)$. Then, see Aubin \cite{Aub},
$$\frac{\int_M\vert\nabla u_{x_0}^\varepsilon\vert^2dv_g 
+ \frac{\Lambda_0(g)}{K_n^2}\int_MA_{ii}(u_{x_0}^\varepsilon)^2dv_g}
{\left(\int_M\vert u_{x_0}^\varepsilon\vert^{2^\star}dv_g\right)^{2/2^\star}} < \frac{1}{K_n^2}$$
if $n \ge 4$, $\varepsilon > 0$ is sufficiently small, and (\ref{2ndCstIneq}) is not satisfied at $x_0$. 
This proves (\ref{2ndCstIneq}). 
Now, in order to end the proof 
of Theorem \ref{TheThm}, it remains to prove the assertions in the theorem 
concerning extremal maps. This is the subject of Sections \ref{SectProof2} and \ref{SectProof4}.

\section{Proof of the second part of Theorem \ref{TheThm}}\label{SectProof2}

As in Section \ref{SectProof1}, we let $(M,g)$ be a smooth compact Riemannian manifold 
of dimension $n \ge 3$, $p \ge 1$ an integer, 
and $A: M \to M_p^s({\mathbb R})$, $A = (A_{ij})$, be 
smooth and such that $A(x)$ is positive as a bilinear form for all $x \in M$. We know from Section 
\ref{SectProof1} that $K_s = K_n^2$, and that (\ref{SharpIneq}) and (\ref{2ndCstIneq})
are true. It remains to prove that if 
the inequality in (\ref{2ndCstIneq}) is strict for all $i$ and all $x$, then the sharp and saturated inequality
(\ref{SharpSatIneq}) 
possesses extremal maps, and the $L^{2^\star}$-normalized set of such extremal maps is precompact in the 
$C^{2,\theta}_p$-topology, where $0 < \theta < 1$. It remains also 
to prove that extremal maps, when they exist, are in general of undeterminate sign, but that 
they can be chosen weakly positive if $-A$ is cooperative, and strongly positive if $A$ is also 
globally irreducible. We claim that the existence of extremal maps when 
the inequality in (\ref{2ndCstIneq}) is strict, and compactness of extremal maps, 
follow from Lemma \ref{TheLem} below. In the sequel, $\Delta_g = -div_g\nabla$ is the Laplace-Beltrami 
operator with respect to $g$.

\begin{lem}\label{TheLem} Let $(M,g)$ be a smooth compact Riemannian manifold 
of dimension $n \ge 4$, $p \ge 1$ an integer, 
and $A_0: M \to M_p^s({\mathbb R})$, $A_0 = (A^0_{ij})$, be  
smooth, such that $A_0(x)$ is positive as a bilinear form for all $x \in M$, and such that 
for any $i$ and any $x$,
$A^0_{ii}(x) > \frac{n-2}{4(n-1)}S_g(x)$, where $S_g$ is the scalar curvature 
of $g$. Let $\left(A(\alpha)\right)_\alpha$, $\alpha\in{\mathbb N}$, be a sequence of smooth maps 
$A(\alpha): M \to M_p^s({\mathbb R})$ such that $A^\alpha_{ij} \to A^0_{ij}$ in $C^{0,\theta}(M)$ as 
$\alpha \to +\infty$, for all $i, j$, where the $A^\alpha_{ij}$'s are the components of $A(\alpha)$, 
and $0 < \theta < 1$. 
Let also $({\mathcal U}_\alpha)_\alpha$  be a sequence of 
$C^{2,\theta}$-solutions of the $p$-systems
\begin{equation}\label{GenericEqtLem}
\Delta_gu_\alpha^i + \sum_{j=1}^pA^\alpha_{ij}(x)u_\alpha^j = 
\lambda_\alpha\vert u_\alpha^i\vert^{2^\star-2}u_\alpha^i
\end{equation}
for all $i$ and all $\alpha$, 
such that $\int_M\vert{\mathcal U}_\alpha\vert^{2^\star}dv_g = 1$ and 
$0 < \lambda_\alpha \le K_n^{-2}$ for all $\alpha$, where the $u_\alpha^i$'s are the components of 
${\mathcal U}_\alpha$. 
Then, up to a subsequence, ${\mathcal U}_\alpha \to {\mathcal U}^0$ 
in $C^{2,\theta}_p(M)$ as $\alpha \to +\infty$, where  ${\mathcal U}^0$ is a nontrivial $p$-map in 
$C^{2,\theta}_p(M)$.
\end{lem}

The proof of Lemma \ref{TheLem} is postponed to Section \ref{SectProof3}. 
We prove here, in Section \ref{SectProof2}, 
that when 
the inequality in (\ref{2ndCstIneq}) is strict for all $i$, 
the existence of extremal maps, and compactness of extremal maps, 
follow from the lemma. The $A(\alpha)$'s in our context are either like 
$A(\alpha) = \Lambda_\alpha K_n^{-2} A$, where the $\Lambda_\alpha$'s are real numbers converging to 
$\Lambda_0(g)$, or like $A(\alpha) = \Lambda_0(g) K_n^{-2}A$ for all $\alpha$, where 
$\Lambda_0(g)$ and $A$ are as in 
Theorem \ref{TheThm}.  Extensions of Lemma \ref{TheLem} 
to higher energies, in the case of conformally flat manifolds, are in Hebey \cite{Heb}. 
The manifold in Lemma \ref{TheLem} needs not to be conformally flat. 
Possible references on elliptic systems, not necessarily like (\ref{GenericEqtLem}), are 
Amster, De N\'apoli, and Mariani \cite{AmsNapMar}, De Figueiredo \cite{DeF}, 
De Figueiredo and Ding \cite{DeFDin}, De Figueiredo and Felmer \cite{DeFFel}, 
Hulshof, Mitidieri and Vandervorst \cite{HulMitVan}, 
Mitidieri and Sweers \cite{MitSwe}, and Sweers \cite{Swe}.

\medskip We assume that Lemma \ref{TheLem} is true. 
Given $\Lambda > 0$, and ${\mathcal U} \in H_{1,p}^2(M)$, we define the energies 
$E_g^\Lambda({\mathcal U})$ and $\Phi_g({\mathcal U})$ by
\begin{equation}\label{DefEner}
E_g^\Lambda({\mathcal U}) = \int_M\vert\nabla{\mathcal U}\vert^2dv_g 
+ \frac{\Lambda}{K_n^2}\int_MA({\mathcal U},{\mathcal U})dv_g
\end{equation}
and $\Phi_g({\mathcal U}) = \int_M\vert{\mathcal U}\vert^{2^\star}dv_g$. 
By definition of $\Lambda_0(g)$,
\begin{equation}\label{Sec2Eqt1}
\inf_{{\mathcal U} \in {\mathcal H}} E_g^\Lambda({\mathcal U}) < \frac{1}{K_n^2}
\end{equation}
when $\Lambda < \Lambda_0(g)$, where ${\mathcal H}$ is the set
consisting of the ${\mathcal U} \in H_{1,p}^2(M)$ which are such that 
$\Phi_g({\mathcal U}) = 1$.
Let $(\Lambda_\alpha)_\alpha$ be a sequence of positive 
real numbers such that $\Lambda_\alpha < \Lambda_0(g)$ for all $\alpha$, and 
$\Lambda_\alpha \to \Lambda_0(g)$ as $\alpha \to +\infty$. Let also 
 $\lambda_\alpha$ be the infimum in (\ref{Sec2Eqt1}) when we let $\Lambda = \Lambda_\alpha$. Since 
$A > 0$ as a bilinear form, $\lambda_\alpha$ is positive for all $\alpha$. 
By the strict inequality in (\ref{Sec2Eqt1}), 
see Hebey \cite{Heb}, for any $\alpha$, 
there exists ${\mathcal U}_\alpha = (u_\alpha^1,\dots,u_\alpha^p)$ a minimizer 
for $\lambda_\alpha$. In particular, the ${\mathcal U}_\alpha$'s are solutions of the $p$-systems
\begin{equation}\label{Sec2Eqt2}
\Delta_gu_\alpha^i + \frac{\Lambda_\alpha}{K_n^2}\sum_{j=1}^pA_{ij}(x)u_\alpha^j = 
\lambda_\alpha\vert u_\alpha^i\vert^{2^\star-2}u_\alpha^i
\end{equation}
for all $i$, and such that $\Phi_g({\mathcal U}_\alpha) = 1$ 
and ${\mathcal U}_\alpha \in C^{2,\theta}_p(M)$ for all $\alpha$, where 
$0 < \theta < 1$. Up to a subsequence, 
we may assume that $\lambda_\alpha \to \lambda_0$ as $\alpha \to +\infty$. If 
the inequality in (\ref{2ndCstIneq}) is strict for all $i$, we can apply Lemma 
\ref{TheLem} with $A(\alpha) = \Lambda_\alpha K_n^{-2} A$, and 
$A_0 = \Lambda_0(g)K_n^{-2}A$. 
By Lemma \ref{TheLem} we then get that, up to a subsequence, the ${\mathcal U}_\alpha$'s converge in $C^{2,\theta}_p(M)$ to 
some ${\mathcal U}^0$. Then $\Phi_g({\mathcal U}^0) = 1$, and, by (\ref{Sec2Eqt2}),
\begin{equation}\label{Sec2Eqt3}
\Delta_gu^0_i + \frac{1}{K_n^2}\sum_{j=1}^pA^0_{ij}(x)u^0_j = 
\lambda_0 \vert u^0_i\vert^{2^\star-2}u^0_i
\end{equation}
for all $i$, where the $A^0_{ij}$'s are the components of the matrix $A_0(g) = \Lambda_0(g)A$. Since 
we have that $\lambda_\alpha < K_n^{-2}$ for all $\alpha$, we can write that $\lambda_0 \le K_n^{-2}$. On the 
other hand, multiplying (\ref{Sec2Eqt3}) by $u^0_i$, integrating over $M$, and summing over $i$, 
we get that
$$E_g^{\Lambda_0(g)}({\mathcal U}^0) = \lambda_0
\hskip.1cm ,$$
where $E^\Lambda_g$ is given by (\ref{DefEner}). 
By the definition of $\Lambda_0(g)$, it follows that $\lambda_0 \ge K_n^{-2}$. In particular, 
$\lambda_0 = K_n^{-2}$, and ${\mathcal U}^0$ is a nontrivial extremal map for 
(\ref{SharpSatIneq}). This proves the above claim that if 
the inequality in (\ref{2ndCstIneq}) is strict for all $i$, 
then the existence of extremal maps follows from Lemma \ref{TheLem}.

\medskip Concerning compactness, let ${\mathcal H}_0$ be the $L^{2^\star}$-normalized set of 
extremal maps for (\ref{SharpSatIneq}). Then ${\mathcal H}_0$ consists of the ${\mathcal U}^0 \in H_{1,p}^2(M)$ 
such that $\Phi_g({\mathcal U}^0) = 1$, and
$$E_g^{\Lambda_0(g)}({\mathcal U}^0) = \inf_{\left\{\Phi_g({\mathcal U}) = 1\right\}} 
E_g^{\Lambda_0(g)}({\mathcal U}) = \frac{1}{K_n^2}
\hskip.1cm ,$$
where $E^\Lambda_g$ is given by (\ref{DefEner}). In particular, 
the extremal maps ${\mathcal U}^0$ in ${\mathcal H}_0$ are solutions of the $p$-system
\begin{equation}\label{Sec2Eqt4}
\Delta_gu^0_i + \frac{1}{K_n^2}\sum_{j=1}^pA^0_{ij}(x)u^0_j = 
K_n^{-2} \vert u^0_i\vert^{2^\star-2}u^0_i
\end{equation}
for all $i$, and such that $\Phi_g({\mathcal U}^0) = 1$, 
where the $A^0_{ij}$'s are the components of the matrix $A_0(g) = \Lambda_0(g)A$, and the $u^0_i$'s  
are the components of ${\mathcal U}$. Such ${\mathcal U}^0$'s, see Hebey \cite{Heb}, are in 
$C^{2,\theta}_p(M)$, where $0 < \theta < 1$. If 
the inequality in (\ref{2ndCstIneq}) is strict for all $i$, we can apply Lemma 
\ref{TheLem} with $A(\alpha) = A_0 = \Lambda_0(g)K_n^{-2}A$. We get  that 
any sequence in ${\mathcal H}_0$ possesses a converging subsequence 
in $C^{2,\theta}_p(M)$. In particular, ${\mathcal H}_0$ is 
precompact in the $C^{2,\theta}_p$-topology. This proves the above claim that 
if the inequality in (\ref{2ndCstIneq}) is strict for all $i$, 
then the compactness of the set of extremal maps in Theorem \ref{TheThm} 
follows from Lemma \ref{TheLem}.

\medskip Now, in order to end this section, we discuss the assertions in Theorem 
\ref{TheThm} concerning the sign of extremal maps. Extremal maps for (\ref{SharpSatIneq}) 
are solutions of systems like
\begin{equation}\label{Sec2Eqt5}
\Delta_gu_i + \sum_{j=1}^pA^0_{ij}(x)u_j = 
\Lambda \vert u_i\vert^{2^\star-2}u_i
\end{equation}
for all $i$, where $A_0 = (A^0_{ij})$ is like $A_0 = tA$ for some $t > 0$, and $\Lambda = K_n^{-2}$ 
is positive. General remarks on weak solutions of (\ref{Sec2Eqt5}) are as follows. First we 
can note (see, for instance, Hebey \cite{Heb}) that 
weak solutions of such systems are in $C^{2,\theta}_p(M)$, $0 < \theta < 1$. 
Then, when $p \ge 2$, and no specific assumption is made on $A_0$, we can note that 
there are no maximum principles for such systems. For instance, 
see again Hebey \cite{Heb}, we can construct examples of $p$-systems like (\ref{Sec2Eqt5}), $p \ge 2$, 
such that the system possesses solutions with the property that the factors 
of the solutions are nonnegative,  nonzero, but with 
zeros in $M$. Such a phenomenon  
does not occur when $p = 1$ since, when $p = 1$, the maximum principle can be applied and nonnegative solutions 
are either identically zero or everywhere positive. On the other hand, we recover 
the maximum principle for (\ref{Sec2Eqt5}) if we assume that $-A_0$ is cooperative. 
Indeed, when $-A_0$ is cooperative, 
nonnegative solutions of (\ref{Sec2Eqt5}) are such that
$$\Delta_gu_i + A^0_{ii}u_i \ge \lambda u_i^{2^\star-1}$$
for all $i$, and the classical maximum principle for functions can be applied so that either $u_i  > 0$ 
everywhere in $M$, or $u_i \equiv 0$. In particular, in this case, nonnegative solutions 
of (\ref{Sec2Eqt5}) are weakly positive. 
Still when $-A_0$ is 
cooperative, if ${\mathcal U}$ is a weakly positive solution 
of the system, with zero factors, then $A_0$ can be factorized in blocs with respect to the zero and 
nonzero components of ${\mathcal U}$. More precisely, if we write ${\mathcal U} = (u_1,\dots,u_k,0,\dots,0)$ 
with $k < p$, and $u_i > 0$ for all $i$, then
\begin{equation}\label{ExSecMatrix}
A_0 = \left(
\begin{matrix}
S & 0\\
0 & T
\end{matrix}
\right)\hskip.1cm ,
\end{equation}
where $S: M \to M_k^s({\mathbb R})$, $T: M \to M_{p-k}^s({\mathbb R})$, and the $0$'s are null matrix of 
respective order $k\times(p-k)$ and $(p-k)\times k$. This easily follows from the equations 
$\sum_{j=1}^kA^0_{ij}u_j = 0$ 
for all $i \ge k+1$, so that we necessarily have that $A^0_{ij} = 0$ for all $i \ge k+1$ 
and $j \le k$. In this case, the $p$-system (\ref{Sec2Eqt5}) splits into two independent systems -- 
a $k$-system where $A_0$ is replaced by $S$, and a $(p-k)$-system 
where $A_0$ is replaced by $T$. In particular, if $-A_0$ is cooperative and 
$A_0$ is globally irreducible, so that (\ref{ExSecMatrix}) cannot be true, 
then any weakly positive solution of the system is also strongly positive.

\medskip Coming back to minimizers, and to Theorem \ref{TheThm}, 
the first assertion concerning the sign 
of extremal maps in Theorem \ref{TheThm} is that extremal maps 
might be of undeterminate sign when no specific assumption is made on $A$. 
Of course this has to be understood when $p \ge 2$ since, 
when $p = 1$, the maximum principle for functions can be applied. When $p = 1$,  
extremal functions  are either positive or negative.
We assume in what follows that $p = 2$, and let $A$, $A^\prime$ be the matrix
\begin{equation}\label{ExPosMinSecMatrix}
A = \left(
\begin{matrix}
\alpha & \beta\\
\beta & \gamma
\end{matrix}
\right)\hskip.2cm\hbox{and}\hskip.2cm
A^\prime = \left(
\begin{matrix}
\alpha & -\beta\\
-\beta & \gamma
\end{matrix}
\right)\hskip.1cm ,
\end{equation}
where $\alpha, \beta, \gamma$ are smooth functions in $M$, and $A(x)$ is positive for all $x$ 
as a bilinear form. 
For ${\mathcal U} = (u,v)$ in $H_{1,2}^2(M)$ 
we let ${\mathcal U}^\prime$ be given by ${\mathcal U}^\prime = (u,-v)$. We let also 
$\beta \ge 0$, $\beta \not\equiv 0$, be such that it is nontrivial and nonnegative. 
Noting that $A({\mathcal U},{\mathcal U}) = A^\prime({\mathcal U}^\prime,{\mathcal U}^\prime)$, we easily get that 
if ${\mathcal U}_0 = (u_0,v_0)$ 
is an extremal map for the sharp and saturated inequality (\ref{SharpSatIneq}), 
then ${\mathcal U}_0^\prime$ is an extremal map for the modified problem we get by replacing 
$A$ by $A^\prime$, where $A$, $A^\prime$ are as in (\ref{ExPosMinSecMatrix}).
Since ${\mathcal U}_0$ is an extremal map for (\ref{SharpSatIneq}), 
it is also a minimizer for 
$F = E/\Phi_g^{2/2^\star}$, where $\Phi_g({\mathcal U}) = \int_M\vert{\mathcal U}\vert^{2^\star}dv_g$, 
$E({\mathcal U}) = E_g^\Lambda({\mathcal U})$ is as in (\ref{DefEner}), and $\Lambda = \Lambda_0(g)$. 
In particular, $F({\mathcal U}_0) \le F({\mathcal U}_0^\prime)$, and it follows that 
\begin{equation}\label{SignEqt}
\int_M\beta u_0v_0dv_g \le 0\hskip.1cm .
\end{equation}
Since $\beta \ge 0$, 
$-A^\prime$ is cooperative, and we can also write that 
$A^\prime(\hat{\mathcal U}_0,\hat{\mathcal U}_0) \le 
A^\prime({\mathcal U}_0^\prime,{\mathcal U}_0^\prime)$, where $\hat{\mathcal U}_0$ is given by  
$\hat{\mathcal U}_0 = 
(\vert u_0\vert,\vert v_0\vert)$. In particular, $\hat{\mathcal U}_0$ is also an extremal map for the 
modified problem we get by replacing $A$ by $A^\prime$. 
Since $\beta\not\equiv 0$, 
$A^\prime$ is globally irreducible, and it follows from the above discussion that 
$\vert u_0\vert$ and $\vert v_0\vert$ are positive functions. Then, 
by (\ref{SignEqt}), ${\mathcal U}_0$ is like ${\mathcal U}_0 = (u_0,-v_0)$ or ${\mathcal U}_0 = (-u_0,v_0)$ 
where $u_0$ and $v_0$ are positive functions. In particular, neither ${\mathcal U}_0$ nor $-{\mathcal U}_0$ 
are nonnegative. Clearly, this type of discussion extends to integers $p \ge 2$. 
For instance, when $p = 3$, choosing $A$ such that $A_{12}, A_{23} \ge 0$ and $A_{13} \le 0$, we easily construct 
minimizers like ${\mathcal U}_0 = (u_0,-v_0,w_0)$ or ${\mathcal U}_0 = (-u_0,v_0,-w_0)$, where $u_0$, $v_0$, $w_0$ 
are positive functions. This proves the above claim that, 
when no specific assumption is made on $A$, extremal maps for (\ref{SharpSatIneq}) 
might be of undeterminate sign. On the contrary, if we assume that $-A$ is cooperative, 
then $A(\hat{\mathcal U},\hat{\mathcal U}) \le 
A({\mathcal U},{\mathcal U})$ for all ${\mathcal U} = (u_1,\dots,u_p)$, where 
$\hat{\mathcal U} = (\vert u_1\vert,\dots,\vert u_p\vert)$. In particular, if 
${\mathcal U}_0$ is an extremal map for (\ref{SharpSatIneq}), then 
$\hat{\mathcal U}_0$ is also an extremal map for (\ref{SharpSatIneq}). 
By the above discussion 
for systems like (\ref{Sec2Eqt5}), $\hat{\mathcal U}_0$ has to be weakly positive since $-A$ is cooperative. It is 
even strongly positive if $A$ is also globally irreducible. 
In particular, 
extremal maps for (\ref{SharpSatIneq}) can be chosen weakly positive 
when $-A$ is cooperative, and even strongly positive 
$A$ is also globally irreducible. This proves the assertions 
in Theorem \ref{TheThm}  concerning the sign 
of extremal maps. Up to Lemma \ref{TheLem}, Theorem \ref{TheThm} is proved.

\medskip When $-A$ is cooperative, and $A$ is globally irreducible, we 
can prove the stronger result that 
any extremal map ${\mathcal U}$ for (\ref{SharpSatIneq}) has to be such 
that either ${\mathcal U}$ or $-{\mathcal U}$ is 
strongly positive. In order to see this we first note that, according to the above 
proof, when $-A$ is cooperative, and $A$ is globally irreducible, 
the components of an 
extremal map for (\ref{SharpSatIneq}) 
are either positive or negative functions. By contradiction, 
up to permuting the indices, we write that ${\mathcal U} = (u_1,\dots,u_k,-u_{k+1},\dots,-u_p)$ 
is an extremal map for (\ref{SharpSatIneq}), where the $u_i$'s are positive functions. We let 
${\mathcal U}^\prime$ be given by ${\mathcal U}^\prime = (u_1,\dots,u_p)$. 
Writing that $E^\Lambda_g({\mathcal U}) \le E^\Lambda_g({\mathcal U}^\prime)$, where 
$E^\Lambda_g$ is as in (\ref{DefEner}) and $\Lambda = \Lambda_0(g)$, 
we get that
$$\sum_{i\in{\mathcal H}_k, j \in {\mathcal H}_{k+1}}  A_{ij}u_iu_j \ge 0\hskip.1cm ,$$
where ${\mathcal H}_k = \left\{1,\dots,k\right\}$, and ${\mathcal H}_{k+1} = \left\{k+1,\dots,p\right\}$. 
The contradiction follows since $-A$ is cooperative, $A$ is globally irreducible, and the $u_i$'s are positive 
functions. This 
proves that when $-A$ is cooperative, and $A$ is globally irreducible, 
extremal maps ${\mathcal U}$ for (\ref{SharpSatIneq}) are such that either ${\mathcal U}$ or $-{\mathcal U}$ is 
strongly positive.

\section{Applications of Theorem \ref{TheThm}}\label{SectProof3}

We discuss the two corollaries, or applications, of Theorem \ref{TheThm} we briefly  mentionned 
in the introduction. The first application, stating that  
the sharp and saturated inequality (\ref{SharpSatIneq}) 
possesses extremal maps when $(M,g)$ has nonpositive scalar curvature and $n \ge 4$, is easy to get. Indeed, 
since $A$ in Theorem \ref{TheThm} 
is such that $A(x)$ is positive in the sense of bilinear forms for all $x$, we clearly have that 
$A_{ii}(x) > 0$ for all $x$ and all $i$. In particular, 
(\ref{2ndCstIneq}) is always true when $(M,g)$ has nonpositive scalar curvature.

\medskip A less obvious result is the second application stating that 
if $n \ge 4$, $A$ does not depend on $x$, and $(M,g)$ has constant scalar curvature, then (\ref{SharpSatIneq}) 
possesses extremal maps. When $(M,g)$ is not conformally diffeomorphic to the unit sphere, the 
result easily follows from the developments in 
Aubin \cite{Aub} and Schoen \cite{Sch}. The energy estimates in 
Aubin \cite{Aub} and Schoen \cite{Sch} give that, in this case, 
when $(M,g)$ is not conformally diffeomorphic to the unit sphere, 
the inequality in (\ref{2ndCstIneq}) has to be strict. Then we can apply Theorem \ref{TheThm}. 
When $(M,g)$ is the unit sphere, or conformally 
diffeomorphic to the unit sphere, the only problem is when 
equality holds in (\ref{2ndCstIneq}) for one, or at least one $i$. For such an $i$, we claim that we necessarily have that  $A_{ij} = 0$ for all $j \not=i$. Assuming for the moment that 
the claim is true, we easily get with such a claim that there exist 
extremal maps for (\ref{SharpSatIneq}). The sharp and saturated scalar Sobolev inequality on the 
unit sphere $(S^n,g_0)$ reads as
$$\left(\int_{S^n}\vert u\vert^{2^\star}dv_{g_0}\right)^{2/2^\star} \le K_n^2\int_{S^n}\vert\nabla u\vert^2dv_{g_0} 
+ \omega_n^{-2/n} \int_{S^n}u^2dv_{g_0}\hskip.1cm ,$$
where $\omega_n$ is the volume of the unit sphere. In particular, see for instance Hebey \cite{HebCIMSBook} for a reference 
in book form, 
there is a whole family of extremal functions for the inequality, including constant functions. Let $u_0$ be one of these 
functions. We choose $u_0$ such that 
$u_0$ is positive and $\Vert u_0\Vert_{2^\star} = 1$. When $(M,g)$ is the unit sphere, equality 
holds in (\ref{2ndCstIneq}) for one $i$, and $A_{ij} = 0$ for all $j \not=i$, the $p$-map ${\mathcal U} = (u_1,\dots,u_p)$, 
where $u_i = u_0$ and $u_j = 0$ for $j \not= i$, is clearly an extremal map for (\ref{SharpSatIneq}). In particular, 
(\ref{SharpSatIneq}) possesses an extremal map. It remains to prove the above claim that when 
$(M,g)$ is the unit sphere, and equality 
holds in (\ref{2ndCstIneq}) for one $i$, we necessarily 
have that $A_{ij} = 0$ for all $j\not= i$. In order to prove this, 
we proceed by contradiction. We assume that $(M,g)$ is the unit sphere, that equality 
holds in (\ref{2ndCstIneq}) for one $i$, and that there exists $j \not= i$ such that $A_{ij} \not= 0$. We let 
${\mathcal U}_\varepsilon = (u_\varepsilon^1,\dots,u_\varepsilon^p)$ be given by 
$u_\varepsilon^i = u_0$, where $u_0$ is as above, 
$u_\varepsilon^j = -\varepsilon A_{ij}$, and $u_\varepsilon^k = 0$ if $k \not= i, j$, 
where $\varepsilon > 0$ is small. Then, with the notations in Theorem \ref{TheThm},
\begin{eqnarray*}
&&K_n^2\int_{S^n}\vert\nabla{\mathcal U}_\varepsilon\vert^2dv_{g_0} 
+ \Lambda_0(g)\int_{S^n}A({\mathcal U}_\varepsilon,{\mathcal U}_\varepsilon)dv_{g_0}\\
&&= K_n^2\int_{S^n}\vert\nabla u_0\vert^2dv_{g_0} + \omega_n^{-2/n} \int_{S^n}u_0^2dv_{g_0}\\
&&\hskip.4cm - 2\Lambda_0(g_0)A_{ij}^2\varepsilon\int_{S^n}u_0dv_{g_0} + O\left(\varepsilon^2\right)\\
&&= 1 - 2\Lambda_0(g_0)A_{ij}^2\varepsilon\int_{S^n}u_0dv_{g_0} + O\left(\varepsilon^2\right)
\end{eqnarray*}
and since $\int_{S^n}\vert{\mathcal U}_\varepsilon\vert^{2^\star}dv_{g_0} \ge 1$, we get a contradiction 
with (\ref{SharpSatIneq}) by choosing  
$\varepsilon > 0$ sufficiently small. This proves the above claim that when 
$(M,g)$ is the unit sphere, and equality 
holds in (\ref{2ndCstIneq}) for one $i$, we cannot have that there exists $j \not= i$ such that $A_{ij} \not= 0$. 
This also ends the proof of the second application of Theorem \ref{TheThm} stating that 
if $n \ge 4$, $A$ does not depend on $x$, and $(M,g)$ has constant scalar curvature, then (\ref{SharpSatIneq}) 
possesses extremal maps.

\section{Proof of Lemma \ref{TheLem}}\label{SectProof4}

We prove Lemma \ref{TheLem} in this Section. 
We let $(M,g)$ be a smooth compact Riemannian manifold 
of dimension $n \ge 3$, $p \ge 1$ an integer, 
and $A_0: M \to M_p^s({\mathbb R})$, $A_0 = (A^0_{ij})$, be  
smooth and such that $A_0(x)$ is positive as a bilinear form for all $x \in M$.
We let $\left(A(\alpha)\right)_\alpha$, $\alpha\in{\mathbb N}$, be a sequence of smooth maps 
$A(\alpha): M \to M_p^s({\mathbb R})$ such that $A^\alpha_{ij} \to A^0_{ij}$ in $C^{0,\theta}(M)$ as 
$\alpha \to +\infty$, for all $i, j$, where the $A^\alpha_{ij}$'s are the components of $A(\alpha)$, 
and $0 < \theta < 1$. 
We let also $({\mathcal U}_\alpha)_\alpha$  be a sequence of 
$C^{2,\theta}$-solutions of the $p$-systems
\begin{equation}\label{GenericEqtLemSec4}
\Delta_gu_\alpha^i + \sum_{j=1}^pA^\alpha_{ij}(x)u_\alpha^j = 
\lambda_\alpha\vert u_\alpha^i\vert^{2^\star-2}u_\alpha^i
\end{equation}
for all $i$ and all $\alpha$, 
such that $\int_M\vert{\mathcal U}_\alpha\vert^{2^\star}dv_g = 1$ and 
$\lambda_\alpha \le K_n^{-2}$ for all $\alpha$, where the $u_\alpha^i$'s are the components of 
${\mathcal U}_\alpha$. Since $A_0(x)$ is positive as a bilinear form for all $x$, and since 
$A^\alpha_{ij} \to A^0_{ij}$ in $C^{0,\theta}(M)$ as $\alpha \to +\infty$ for all $i$, $j$, there exists 
$K > 0$ such that $A^\alpha_{ij}(x) \ge K\delta_{ij}$ in the sense of bilinear forms, for all $x$ and 
all $\alpha$ sufficiently large. Multiplying (\ref{GenericEqtLemSec4}) by $u_\alpha^i$, integrating over $M$, and summing 
over $i$, we then get with the Sobolev inequality that there exists $\lambda > 0$ such that 
$\lambda_\alpha \ge \lambda$ for all $\alpha$ sufficiently large. Now we define $\tilde{\mathcal U}_\alpha$ by 
$\tilde{\mathcal U}_\alpha = (\tilde u_\alpha^1,\dots,\tilde u_\alpha^p)$, where
\begin{equation}\label{Sec4ProofEqt1}
\tilde u_\alpha^i = \lambda_\alpha^{\frac{n-2}{4}}u_\alpha^i
\end{equation}
for all $\alpha$ and all $i$. Then, for any $\alpha$, 
$\tilde{\mathcal U}_\alpha$ is a solution of the $p$-system
\begin{equation}\label{Sec4ProofEqt2}
\Delta_g\tilde u_\alpha^i + \sum_{j=1}^pA^\alpha_{ij}(x)\tilde u_\alpha^j = 
\vert\tilde u_\alpha^i\vert^{2^\star-2}\tilde u_\alpha^i
\end{equation}
for all $i$. Moreover, since $\lambda_\alpha \le K_n^{-2}$, we also have that
\begin{equation}\label{Sec4ProofEqt3}
\int_M\vert\tilde{\mathcal U}_\alpha\vert^{2^\star}dv_g \le K_n^{-n}
\end{equation}
for all $\alpha$. Lemma \ref{TheLem} states that, up to a subsequence, ${\mathcal U}_\alpha \to {\mathcal U}^0$ 
in $C^{2,\theta}_p(M)$ as $\alpha \to +\infty$, where  ${\mathcal U}^0$ is a nontrivial $p$-map in 
$C^{2,\theta}_p(M)$. By standard elliptic theory, and (\ref{GenericEqtLemSec4}), it suffices to prove that, up to a subsequence,  
the ${\mathcal U}_\alpha$'s are uniformly bounded in $L^\infty(M)$. Since $\lambda_\alpha \in [\lambda,K_n^{-2}]$ 
for $\alpha$ large, the  ${\mathcal U}_\alpha$'s are uniformly bounded in $L^\infty(M)$ if and only if the 
$\tilde{\mathcal U}_\alpha$'s are uniformly bounded in $L^\infty(M)$. In particular, Lemma \ref{TheLem} reduces to 
proving that, up to a subsequence, there exists $C > 0$ such that
\begin{equation}\label{MainResToProve}
\vert\tilde{\mathcal U}_\alpha\vert \le C
\end{equation}
in $M$, for all $\alpha$, where $\vert\tilde{\mathcal U}_\alpha\vert = \sum_i\vert\tilde u_\alpha^i\vert$. We prove 
(\ref{MainResToProve}) in what follows. Up to a subsequence, by compactness of the embedding 
$H_1^2 \subset L^2$, we may assume that 
$\tilde{\mathcal U}_\alpha \to \tilde{\mathcal U}^0$ in $L^2$ for some $\tilde{\mathcal U}^0 \in H_{1,p}^2(M)$. 
In other words, we can assume that there are functions $\tilde u^0_i$ in $H_1^2(M)$ such that
\begin{equation}\label{Sec4ProofEqt4}
\tilde u_\alpha^i \to \tilde u^0_i\hskip.2cm\hbox{in}\hskip.1cm L^2(M)
\end{equation}
for all $i$, as $\alpha \to +\infty$. We may also assume that $\tilde u_\alpha^i \rightharpoonup \tilde u^0_i$ 
weakly in $H_1^2(M)$, that $\tilde u_\alpha^i \to \tilde u^0_i$ a.e in $M$, and that 
$\vert\tilde u_\alpha^i\vert^{2^\star-2}\tilde u_\alpha^i 
\rightharpoonup \vert\tilde u^0_i\vert^{2^\star-2}\tilde u^0_i$ weakly in $L^{2^\star/(2^\star-1)}(M)$ for all 
$i$. In particular, $\tilde{\mathcal U}^0$ is a solution of the limit equation
\begin{equation}\label{Sec4ProofEqt5}
\Delta_g\tilde u^0_i + \sum_{j=1}^pA^\alpha_{ij}(x)\tilde u^0_j = 
\vert\tilde u^0_i\vert^{2^\star-2}\tilde u^0_i
\end{equation}
for all $i$. Then, see, for instance, Hebey \cite{Heb}, we can prove that $\tilde{\mathcal U}^0$ is in 
$C^{2,\theta}_p(M)$.
For $(x_\alpha)_\alpha$ a converging sequence 
of points in $M$, and $(\mu_\alpha)_\alpha$ a sequence of positive real numbers converging to zero, we define a 
{\it $1$-bubble} as a sequence $(B_\alpha)_\alpha$ of functions in $M$ given by
\begin{equation}\label{Def1BubbleSec4}
B_\alpha(x) = \left(\frac{\mu_\alpha}{\mu_\alpha^2 + \frac{d_g(x_\alpha,x)^2}{n(n-2)}}\right)^{\frac{n-2}{2}}\hskip.1cm .
\end{equation}
The $x_\alpha$'s are referred to as the {\it centers} and the $\mu_\alpha$'s as 
the {\it weights} of the $1$-bubble $(B_\alpha)_\alpha$. We define a {\it $p$-bubble} as a 
sequence $({\mathcal B}_\alpha)_\alpha$ 
of $p$-maps such that, if we write that ${\mathcal B}_\alpha = (B_\alpha^1,\dots,B_\alpha^p)$, then 
$(B_\alpha^i)_\alpha$ is a $1$-bubble for exactly one $i$, and for $j \not= i$, 
$(B_\alpha^j)_\alpha$ is the trivial zero sequence. In other words, a $p$-bubble is a sequence 
of $p$-maps such that one of the components of the sequence is a $1$-bubble, and the other 
components are trivial zero sequences. 
One remark with respect to the definition 
(\ref{Def1BubbleSec4}) is that if $u: {\mathbb R}^n \to {\mathbb R}$ is given by
\begin{equation}\label{PosSolCritEuclEqt}
u(x) = \left(1 + \frac{\vert x\vert^2}{n(n-2)}\right)^{-\frac{n-2}{2}}\hskip.1cm ,
\end{equation}
then $u$ is a positive solution of the critical Euclidean equation $\Delta u = u^{2^\star-1}$, where 
$\Delta = -\sum\partial^2/\partial x_i^2$. More precisely, $u$ is 
the only positive solution of the equation in ${\mathbb R}^n$ 
which is such that $u(0) = 1$ and $u$ is maximum at $0$. All the 
other positive solutions of the equation $\Delta u = u^{2^\star-1}$ in ${\mathbb R}^n$, see 
Caffarelli, Gidas and Spruck \cite{CafGidSpr} and 
Obata \cite{Oba}, 
are then given by
$$\tilde u(x) = \lambda^{(n-2)/2}u\left(\lambda(x-a)\right)\hskip.1cm ,$$
where $\lambda > 0$, and $a \in {\mathbb R}^n$. We prove (\ref{MainResToProve}), and thus 
Lemma \ref{TheLem}, in several steps. 
The first step in the proof is as follows.

\begin{Step}\label{Step1ProofLemSec4} Let $\tilde{\mathcal U}_\alpha$ and $\tilde{\mathcal U}^0$ be 
given by (\ref{Sec4ProofEqt1}) and (\ref{Sec4ProofEqt4}). If $\tilde{\mathcal U}^0 \not\equiv 0$, then 
(\ref{MainResToProve}) is true. If, on the contrary, $\tilde{\mathcal U}^0 \equiv 0$, then 
there exists a $p$-bubble $({\mathcal B}_\alpha)_\alpha$ such that, up to a subsequence, 
\begin{equation}\label{Step1StatemSec4Eqt1}
\tilde{\mathcal U}_\alpha = {\mathcal B}_\alpha + {\mathcal R}_\alpha
\end{equation}
for all $\alpha$, where ${\mathcal R}_\alpha \to 0$ strongly in $H_{1,p}^2(M)$ as $\alpha \to +\infty$. There also exists 
$C > 0$ such that
\begin{equation}\label{Step1StatemSec4Eqt2}
d_g(x_\alpha,x)^{\frac{n-2}{2}}\sum_{i=1}^p\vert\tilde u^i_\alpha(x)\vert \le C
\end{equation}
for all $\alpha$ and all $x \in M$, where the $x_\alpha$'s are the centers of the $1$-bubble from which 
the $p$-bubble $({\mathcal B}_\alpha)_\alpha$ is defined. In particular, 
the $\vert\tilde{\mathcal U}_\alpha\vert$'s are uniformly bounded in any compact subset 
of $M\backslash\{x_0\}$, and $\tilde u^i_\alpha \to 0$ 
in $C^0_{loc}(M\backslash\{x_0\})$ for all $i$ as $\alpha \to +\infty$, where $x_0$ is the limit of the 
$x_\alpha$'s.
\end{Step}

\begin{proof}[Proof of Step \ref{Step1ProofLemSec4}] By the $H_1^2$-theory for blow-up, see Hebey \cite{Heb}, 
there are generalized $p$-bubbles $(\hat{\mathcal B}_{j,\alpha})_\alpha$, 
$j = 1,\dots,k$, such that, up to a subsequence, 
\begin{equation}\label{ProofStep1Sec4Eqt1}
\tilde{\mathcal U}_\alpha = \tilde{\mathcal U}^0 
+ \sum_{j=1}^k\hat{\mathcal B}_{j,\alpha} + {\mathcal R}_\alpha
\end{equation}
and such that
\begin{equation}\label{ProofStep1Sec4Eqt2}
\frac{1}{n}\int_M\vert\tilde{\mathcal U}_\alpha\vert^{2^\star}dv_g = 
\frac{1}{n}\int_M\vert\tilde{\mathcal U}^0\vert^{2^\star}dv_g 
+ \sum_{j=1}^kE_f(\hat{\mathcal B}_{j,\alpha}) + o(1)
\end{equation}
for all $\alpha$, where ${\mathcal R}_\alpha \to 0$ strongly in $H_{1,p}^2(M)$ as $\alpha \to +\infty$, 
$E_f(\hat{\mathcal B}_{j,\alpha})$ is the energy of the generalized $p$-bubble $(\hat {\mathcal B}_{j,\alpha})_\alpha$, 
and $o(1) \to 0$ as $\alpha \to +\infty$. Generalized $p$-bubbles are rescaling of solutions of the critical 
equation $\Delta u = \vert u\vert^{2^\star-2}u$, and the energy of the 
generalized $p$-bubble is the energy of $u$. In particular, 
the energy $E_f(\hat{\mathcal B}_{j,\alpha})$ does not depend on 
$\alpha$. It is 
always greater than or equal to $K_n^{-n}/n$, and if equality holds, then, up 
to lower order terms, the generalized $p$-bubble has to be  
a $p$-bubble. Namely,  we always have that 
$E_f(\hat{\mathcal B}_{j,\alpha}) \ge K_n^{-n}/n$, and if equality holds, then
$$\hat{\mathcal B}_{j,\alpha} = {\mathcal B}_{j,\alpha} + {\mathcal R}_\alpha\hskip.1cm ,$$
where 
$({\mathcal B}_{j,\alpha})_\alpha$ is a $p$-bubble, as defined above, and 
${\mathcal R}_\alpha \to 0$ strongly in $H_{1,p}^2(M)$ as $\alpha \to +\infty$. By 
(\ref{Sec4ProofEqt3}), it follows from (\ref{ProofStep1Sec4Eqt2}) that if $\tilde{\mathcal U}^0 \not\equiv 0$,
then $k = 0$, and that if $\tilde{\mathcal U}^0 \equiv 0$, then $k = 0$ or $k = 1$. When $\tilde{\mathcal U}^0 \equiv 0$, 
and $k = 0$, we get with (\ref{ProofStep1Sec4Eqt1}) that $\tilde{\mathcal U}_\alpha\to 0$ 
strongly in $H_{1,p}^2(M)$ as $\alpha \to +\infty$. This is impossible since, by construction of the 
$\tilde{\mathcal U}_\alpha$'s, we also have that there is a uniform positive lower bound for 
the left hand side in (\ref{ProofStep1Sec4Eqt2}). In particular, $k = 1$ 
when $\tilde{\mathcal U}^0 \equiv 0$. When $\tilde{\mathcal U}^0 \equiv 0$, 
and $k = 1$, we also get from (\ref{Sec4ProofEqt3}), 
(\ref{ProofStep1Sec4Eqt2}), and the above discussion, that the 
generalized $p$-bubble in (\ref{ProofStep1Sec4Eqt1}) has to be a $p$-bubble, 
and that $\lambda_\alpha \to K_n^{-2}$ as $\alpha \to +\infty$.
Summarizing, we get with the  $H_1^2$-theory for blow-up that if $\tilde{\mathcal U}^0 \not\equiv 0$, then, up to a subsequence, 
\begin{equation}\label{ProofStep1Sec4Eqt3}
\tilde{\mathcal U}_\alpha= \tilde{\mathcal U}^0 + {\mathcal R}_\alpha
\end{equation}
for all $\alpha$, and that if $\tilde{\mathcal U}^0 \equiv 0$, 
then 
there exists a $p$-bubble $({\mathcal B}_\alpha)_\alpha$ such that, up to a subsequence,
\begin{equation}\label{ProofStep1Sec4Eqt4}
\tilde{\mathcal U}_\alpha = {\mathcal B}_\alpha + {\mathcal R}_\alpha
\end{equation}
for all $\alpha$, where, in (\ref{ProofStep1Sec4Eqt3}) and (\ref{ProofStep1Sec4Eqt4}), 
${\mathcal R}_\alpha \to 0$ strongly in $H_{1,p}^2(M)$ as $\alpha \to +\infty$. We let 
the $x_\alpha$'s and $\mu_\alpha$'s be the centers and weights of the $1$-bubble from which 
the $p$-bubble $({\mathcal B}_\alpha)_\alpha$ 
in (\ref{ProofStep1Sec4Eqt4}) is defined. We claim that
\begin{equation}\label{ProofStep1Sec4Eqt5}
\hbox{(\ref{MainResToProve}) is true if}\hskip.1cm \tilde{\mathcal U}^0 \not\equiv 0
\hskip.1cm ,\hskip.1cm\hbox{while (\ref{Step1StatemSec4Eqt2}) is true if}\hskip.1cm \tilde{\mathcal U}^0 \equiv 0
\hskip.1cm .
\end{equation}
In order to prove (\ref{ProofStep1Sec4Eqt5}), we let $\Phi_\alpha$ be the function given by 
$\Phi_\alpha(x) = 1$ if $\tilde{\mathcal U}^0 \not\equiv 0$, and 
$\Phi_\alpha(x) = d_g(x_\alpha,x)$ if $\tilde{\mathcal U}^0 \equiv 0$. We let also $\Psi_\alpha$ be the function 
given by
\begin{equation}\label{PsiDefProofThPointEst}
\Psi_\alpha(x) = 
\Phi_\alpha(x)^{\frac{n-2}{2}}
\sum_{i=1}^p\vert\tilde u_\alpha^i(x)\vert
\hskip.1cm .
\end{equation}
Then (\ref{ProofStep1Sec4Eqt5}) is equivalent to the statement that the 
$\Psi_\alpha$'s are uniformly bounded in $L^\infty(M)$. Now we proceed by contradiction. 
We let the $y_\alpha$'s be points in $M$ such that 
the $\Psi_\alpha$'s are maximum at $y_\alpha$ and $\Psi_\alpha(y_\alpha) \to +\infty$ as $\alpha \to +\infty$. 
Up to a subsequence, we may assume that $\vert\tilde u_\alpha^{i_0}\vert(y_\alpha) \ge 
\vert\tilde u_\alpha^i\vert(y_\alpha)$ 
for some $i_0=1,\dots,p$, and all $i$. We set 
$\tilde\mu_\alpha = \vert\tilde u_\alpha^{i_0}\vert(y_\alpha)^{-2/(n-2)}$. Then 
$\tilde\mu_\alpha \to 0$ as $\alpha \to +\infty$, and by (\ref{PsiDefProofThPointEst}) we also have that
\begin{equation}\label{PsiDefProofThPointEstEqt1}
\frac{d_g(x_\alpha,y_\alpha)}{\tilde\mu_\alpha} \to +\infty
\end{equation}
if $\tilde{\mathcal U}^0 \equiv 0$, 
as $\alpha \to +\infty$. Let $\delta > 0$ be less than the injectivity radius of 
$(M,g)$. For $i = 1,\dots,p$, we define the function $\tilde v_\alpha^i$ in $B_0(\delta\mu_\alpha^{-1})$ by
\begin{equation}\label{PsiDefProofThPointEstEqt2}
\tilde v_\alpha^i(x) = 
\tilde\mu_\alpha^{\frac{n-2}{2}}\tilde u_\alpha^i\left(\exp_{y_\alpha}(\tilde\mu_\alpha x)\right)
\hskip.1cm ,
\end{equation}
where $B_0(\delta\mu_\alpha^{-1})$ is the Euclidean ball of radius $\delta\tilde\mu_\alpha^{-1}$ 
centered at $0$, and $\exp_{y_\alpha}$ is the exponential map at $y_\alpha$. 
Given $R > 0$ 
and $x \in B_0(R)$, the Euclidean ball of radius $R$ centered at $0$, we can write with 
(\ref{PsiDefProofThPointEst}) and (\ref{PsiDefProofThPointEstEqt2}) that
\begin{equation}\label{PsiDefProofThPointEstEqt3}
\vert\tilde v_\alpha^i\vert(x) \le \frac{\tilde\mu_\alpha^{\frac{n-2}{2}}\Psi_\alpha\left(\exp_{y_\alpha}(\tilde\mu_\alpha x)\right)}
{\Phi_\alpha\left(\exp_{y_\alpha}(\tilde\mu_\alpha x)\right)^{\frac{n-2}{2}}}
\end{equation}
for all $i$, when $\alpha$ is sufficiently large. 
For any $x \in B_0(R)$, when $\tilde{\mathcal U}^0 \equiv 0$, 
\begin{eqnarray*} d_g\left(x_\alpha,\exp_{y_\alpha}(\tilde \mu_\alpha x)\right)
& \ge & d_g\left(x_\alpha,y_\alpha\right) - R\tilde\mu_\alpha\\
& \ge & \left(1 - \frac{R\tilde\mu_\alpha}{\Phi_\alpha(y_\alpha)}\right)\Phi_\alpha(y_\alpha)
\end{eqnarray*}
when $\alpha$ is sufficiently large so that, by 
(\ref{PsiDefProofThPointEstEqt1}), the right hand side of the last equation is positive. 
Coming back to 
(\ref{PsiDefProofThPointEstEqt3}), thanks to the definition of the $y_\alpha$'s, we then get that 
for any $i$, and any $x \in B_0(R)$,
\begin{equation}\label{PsiDefProofThPointEstEqt4}
\begin{split}
\vert\tilde v_\alpha^i\vert(x)
& \le \frac{\tilde\mu_\alpha^{\frac{n-2}{2}}\Psi_\alpha(y_\alpha)}
{\Phi_\alpha\left(\exp_{y_\alpha}(\tilde\mu_\alpha x)\right)^{\frac{n-2}{2}}}\\
& \le p \left(1 - \frac{R\tilde\mu_\alpha}{\Phi_\alpha(y_\alpha)}\right)^{-\frac{n-2}{2}}
\end{split}
\end{equation}
when $\alpha$ is sufficiently large. 
In particular, by (\ref{PsiDefProofThPointEstEqt1}) and 
(\ref{PsiDefProofThPointEstEqt4}), up to passing 
to a subsequence, the $\tilde v_\alpha^i$'s are uniformly bounded in any compact subset 
of ${\mathbb R}^n$ for all $i$. 
Let $\tilde{\mathcal V}_\alpha = (\tilde v_\alpha^1,\dots,\tilde v_\alpha^p)$. The 
$\tilde{\mathcal V}_\alpha$'s are solutions of the system
\begin{equation}\label{PsiDefProofThPointEstEqt5}
\Delta_{g_\alpha}\tilde v_\alpha^i + \sum_{j=1}^p\tilde\mu_\alpha^2\tilde A_{ij}^\alpha\tilde v_\alpha^j 
= \vert\tilde v_\alpha^i\vert^{2^\star-2}\tilde v_\alpha^i
\hskip.1cm ,
\end{equation}
for all $i$, where
\begin{eqnarray*}
&&\tilde A_{ij}^\alpha(x) = A_{ij}^\alpha\left(\exp_{y_\alpha}(\tilde\mu_\alpha x)\right)
\hskip.1cm ,\hskip.1cm\hbox{and}\\
&&g_\alpha(x) 
= \left(\exp_{y_\alpha}^\star g\right)(\tilde\mu_\alpha x)\hskip.1cm .
\end{eqnarray*}
Let $\xi$ be the Euclidean metric. Clearly, for any compact subset 
$K$ of ${\mathbb R}^n$, $g_\alpha \to \xi$ in $C^2(K)$ as $\alpha \to +\infty$. Then, by standard elliptic theory, 
and (\ref{PsiDefProofThPointEstEqt5}), 
we get that the $\tilde v_\alpha^i$'s are uniformly bounded 
in $C^{2,\theta}_{loc}({\mathbb R}^n)$ for all $i$, where $0 < \theta < 1$. In particular, up to a subsequence, 
we can assume that $\tilde v_\alpha^i \to \tilde v_i$ in $C^2_{loc}({\mathbb R}^n)$ as $\alpha \to +\infty$ 
for all $i$, where the $\tilde v_i$'s 
are functions in $C^2({\mathbb R}^n)$. The $\tilde v_i$'s are bounded in ${\mathbb R}^n$ by 
(\ref{PsiDefProofThPointEstEqt4}), and such that $\vert\tilde v_{i_0}\vert(0) = 1$ by construction. Without loss of generality, we 
may also assume that the $\tilde v_i$'s are in ${\mathcal D}_1^2({\mathbb R}^n)$ and in $L^{2^\star}({\mathbb R}^n)$ 
for all $i$, where 
${\mathcal D}_1^2({\mathbb R}^n)$ is the Beppo-Levi space defined as the completion of $C^\infty_0({\mathbb R}^n)$, 
the space of smooth functions with compact support in ${\mathbb R}^n$, with respect to the norm 
$\Vert u\Vert = \Vert\nabla u\Vert_2$. We let $\tilde{\mathcal V} = (\tilde v_1,\dots,\tilde v_p)$. According to the above, 
$\tilde {\mathcal V} \not\equiv 0$. 
By construction, for any $R > 0$, 
$$\int_{B_{y_\alpha}(R\tilde \mu_\alpha)}\vert\tilde{\mathcal U}_\alpha\vert^{2^\star}dv_g 
= \int_{B_0(R)}\vert\tilde{\mathcal V}_\alpha\vert^{2^\star}dv_{g_\alpha}\hskip.1cm .$$
It follows that for any $R > 0$,
\begin{equation}\label{PsiDefProofThPointEstEqt6}
\int_{B_{y_\alpha}(R\tilde\mu_\alpha)}\vert\tilde{\mathcal U}_\alpha\vert^{2^\star}dv_g 
= \int_{{\mathbb R}^n}\vert\tilde{\mathcal V}\vert^{2^\star}dx + \varepsilon_R(\alpha)\hskip.1cm ,
\end{equation}
where $\varepsilon_R(\alpha)$ is such that $\lim_R\lim_\alpha\varepsilon_R(\alpha) = 0$, 
and the limits are as $\alpha\to+\infty$ and $R\to+\infty$. 
When $\tilde{\mathcal U}^0 \equiv 0$, see for instance Hebey \cite{Heb}, 
we also get with (\ref{PsiDefProofThPointEstEqt1}) that
\begin{equation}\label{PsiDefProofThPointEstEqt7}
\lim_{\alpha\to+\infty}\int_{B_{y_\alpha}(R\tilde\mu_\alpha)}
\vert\tilde{\mathcal B}_\alpha\vert^{2^\star}dv_g = 0
\end{equation}
for all $R > 0$, where $(\tilde{\mathcal B}_\alpha)_\alpha$ is the $p$-bubble in (\ref{ProofStep1Sec4Eqt4}). 
By (\ref{ProofStep1Sec4Eqt3}) and (\ref{ProofStep1Sec4Eqt4}), 
\begin{equation}\label{PsiDefProofThPointEstEqt8}
\begin{split}
\int_{B_{y_\alpha}(R\tilde\mu_\alpha)}\vert\tilde{\mathcal U}_\alpha\vert^{2^\star}dv_g
&\le C\int_{B_{y_\alpha}(R\tilde\mu_\alpha)}\vert\tilde{\mathcal U}^0\vert^{2^\star}dv_g + o(1)\\
&= o(1)
\end{split}
\end{equation}
for all $\alpha$ and $R > 0$ if $\tilde{\mathcal U}^0 \not\equiv 0$, while
\begin{equation}\label{PsiDefProofThPointEstEqt9}
\int_{B_{y_\alpha}(R\tilde\mu_\alpha)}\vert\tilde{\mathcal U}_\alpha\vert^{2^\star}dv_g
\le C\int_{B_{y_\alpha}(R\tilde\mu_\alpha)}\vert\tilde{\mathcal B}_\alpha\vert^{2^\star}dv_g + o(1)
\end{equation}
for all $\alpha$ and $R > 0$ if $\tilde{\mathcal U}^0 \equiv 0$, where 
$C > 0$ is independent of $\alpha$ and $R$, and $o(1) \to 0$ as $\alpha \to +\infty$. 
Combining (\ref{PsiDefProofThPointEstEqt6})--(\ref{PsiDefProofThPointEstEqt9}), letting $\alpha \to +\infty$, 
and then $R\to+\infty$, we get that
$$\int_{{\mathbb R}^n}\vert\tilde{\mathcal V}\vert^{2^\star}dx = 0\hskip.1cm ,$$
and this is in contradiction with the equation $\tilde {\mathcal V} \not\equiv 0$. In particular, 
the $\Psi_\alpha$'s are uniformly bounded in $L^\infty(M)$. This proves (\ref{ProofStep1Sec4Eqt5}). 
When
$\tilde{\mathcal U}^0 \equiv 0$, (\ref{ProofStep1Sec4Eqt5}) gives that (\ref{Step1StatemSec4Eqt2}) 
is true, and if $x_0$ is the limit of the $x_\alpha$'s, (\ref{Step1StatemSec4Eqt2}) gives that the 
$\vert\tilde{\mathcal U}_\alpha\vert$'s are uniformly bounded in any compact subset of $M\backslash\{x_0\}$. 
By standard elliptic theory, equation (\ref{Sec4ProofEqt2}) satisfied by the 
$\tilde{\mathcal U}_\alpha$'s, and since 
$\tilde{\mathcal U}_\alpha \to 0$ in $L^2$ when $\tilde{\mathcal U}^0 \equiv 0$, we get that 
$\vert\tilde{\mathcal U}_\alpha\vert \to 0$ in $C^0_{loc}(M\backslash\{x_0\})$ as 
$\alpha \to +\infty$. This ends the proof of Step \ref{Step1ProofLemSec4}.
 \end{proof}

According to Step \ref{Step1ProofLemSec4}, in order to prove (\ref{MainResToProve}), 
it suffices to prove that the $p$-map 
$\tilde{\mathcal U}^0$ given by (\ref{Sec4ProofEqt4}) 
is not identically zero. We proceed here by contradiction and assume that 
$\tilde{\mathcal U}^0 \equiv 0$. The next step in the proof of (\ref{MainResToProve}) 
consists in proving that the $\tilde{\mathcal U}_\alpha$'s satisfy 
perturbed De Giorgi-Nash-Moser type estimates. 
Step \ref{Step2ProofLemSec4} in the proof of (\ref{MainResToProve}) is as follows. 

\begin{Step}\label{Step2ProofLemSec4} Let $\tilde{\mathcal U}_\alpha$ and $\tilde{\mathcal U}^0$ be 
given by (\ref{Sec4ProofEqt1}) and (\ref{Sec4ProofEqt4}). Assume $\tilde{\mathcal U}^0 \equiv 0$. 
 For any $\delta > 0$, there exists $C > 0$ such that, up to a subsequence,
 \begin{equation}\label{Step2StatemSec4Eqt1}
 \max_{M\backslash B_\delta}\vert\tilde{\mathcal U}_\alpha\vert \le C 
 \int_M\left(1 + \vert\tilde{\mathcal U}_\alpha\vert^{2^\star-2}\right)\vert\tilde{\mathcal U}_\alpha\vert dv_g 
 \end{equation}
for all $\alpha$, where $B_\delta = B_{x_0}(\delta)$ is the ball centered at $x_0$ of 
radius $\delta$, $\vert\tilde{\mathcal U}_\alpha\vert = \sum_{i=1}^p\vert\tilde u_\alpha^i\vert$, 
$\vert\tilde{\mathcal U}_\alpha\vert^{2^\star-2} = \sum_{i=1}^p\vert\tilde u_\alpha^i\vert^{2^\star-2}$, 
and $x_0$ is the limit of the centers of the $1$-bubble from which 
the $p$-bubble $({\mathcal B}_\alpha)_\alpha$ in (\ref{Step1StatemSec4Eqt1}) is defined.
\end{Step}

\begin{proof}[Proof of Step \ref{Step2ProofLemSec4}] Let $B = B_x(r)$ be such that 
$B_x(2r) \subset M\backslash\{x_0\}$. 
By (\ref{Sec4ProofEqt2}), and Step \ref{Step1ProofLemSec4}, $\vert\Delta_g\tilde u_\alpha^i\vert 
\le C\vert\tilde{\mathcal U}_\alpha\vert$ in $B$, for all $i$ and all $\alpha$, where 
$C > 0$ is independent 
of $\alpha$ and $i$. Then we also have that
\begin{equation}\label{Step2Sec4Eqt1}
\left\vert\Delta_g\tilde u_\alpha^i + \tilde u_\alpha^i\right\vert 
\le C\vert\tilde{\mathcal U}_\alpha\vert
\end{equation}
in $B$, for all $i$ and all $\alpha$, where $C > 0$ is independent 
of $\alpha$ and $i$. We define the $\hat u_\alpha^i$'s by
\begin{equation}\label{Step2Sec4Eqt2}
\Delta_g\hat u_\alpha^i + \hat u_\alpha^i = 
\left\vert\Delta_g\tilde u_\alpha^i + \tilde u_\alpha^i\right\vert
\end{equation}
in $M$, for all $\alpha$ and all $i$. Since
$$\Delta_g\left(\hat u_\alpha^i \pm \tilde u_\alpha^i\right) 
+ \left(\hat u_\alpha^i \pm \tilde u_\alpha^i\right) \ge 0\hskip.1cm ,$$
we can write that $\hat u_\alpha^i \ge \vert\tilde u_\alpha^i\vert$ in $M$, for all $\alpha$ and all $i$. 
In particular, the $\hat u_\alpha^i$'s are nonnegative. 
By (\ref{Step2Sec4Eqt1}) and (\ref{Step2Sec4Eqt2}) we also have that
\begin{equation}\label{Step2Sec4Eqt3}
\Delta_g\vert\hat{\mathcal U}_\alpha\vert \le C \vert\hat{\mathcal U}_\alpha\vert
\end{equation}
in $B$, for all $\alpha$, 
where $\hat{\mathcal U}_\alpha$ is the $p$-map of components the $\hat u_\alpha^i$'s, 
$\vert\hat{\mathcal U}_\alpha\vert = \sum_{i=1}^p\hat u_\alpha^i$ 
since the $\hat u_\alpha^i$'s are nonnegative, and $C > 0$ is independent of $\alpha$. 
It easily follows from (\ref{Step2Sec4Eqt1}), (\ref{Step2Sec4Eqt2}), and 
Step \ref{Step1ProofLemSec4} that the $\vert\hat{\mathcal U}_\alpha\vert$'s are uniformly bounded 
in $L^\infty(B)$. Then, by (\ref{Step2Sec4Eqt3}), we can apply the De Giorgi-Nash-Moser iterative scheme 
for functions to the $\vert\hat{\mathcal U}_\alpha\vert$'s. In particular, 
we can write that
\begin{equation}\label{Step2Sec4Eqt4}
\max_{B_x(r/4)}\vert\hat{\mathcal U}_\alpha\vert \le C \int_{B_x(r/2)}\vert\hat{\mathcal U}_\alpha\vert dv_g
\hskip.1cm ,
\end{equation}
where $C > 0$ is independent of $\alpha$. With the notations 
in the statement of Step \ref{Step2ProofLemSec4}, since $B$ is basically any ball 
in $M\backslash\{x_0\}$, it easily follows form (\ref{Step2Sec4Eqt4}) that 
for any $\delta > 0$, 
\begin{equation}\label{Step2Sec4Eqt5}
\max_{M\backslash B_\delta}\vert\hat{\mathcal U}_\alpha\vert \le C 
\int_M\vert\hat{\mathcal U}_\alpha\vert dv_g
\hskip.1cm ,
\end{equation}
where $C > 0$ is independent of $\alpha$. By 
(\ref{Sec4ProofEqt2}), 
\begin{equation}\label{Step2Sec4Eqt6}
\left\vert\Delta_g\tilde u_\alpha^i + \tilde u_\alpha^i\right\vert 
\le C\left(1 + \vert\tilde{\mathcal U}_\alpha\vert^{2^\star-2}\right)\vert\tilde{\mathcal U}_\alpha\vert
\end{equation}
in $M$, for all $i$ and $\alpha$, where $C > 0$ is independent of $\alpha$ and $i$.
Integrating (\ref{Step2Sec4Eqt2}) over $M$, since $\int_M(\Delta_g\hat u_\alpha^i)dv_g = 0$ 
for all $i$ and all $\alpha$, we get with 
(\ref{Step2Sec4Eqt6}) that
\begin{equation}\label{Step2Sec4Eqt7}
\int_M\vert\hat{\mathcal U}_\alpha\vert dv_g \le C
\int_M\left(1 + \vert\tilde{\mathcal U}_\alpha\vert^{2^\star-2}\right)\vert\tilde{\mathcal U}_\alpha\vert dv_g
\end{equation}
for all $\alpha$, where $C > 0$ is independent of $\alpha$. As already mentionned, $\vert\tilde{\mathcal U}_\alpha\vert 
\le \vert\hat{\mathcal U}_\alpha\vert$ in $M$. 
In particular, we get with (\ref{Step2Sec4Eqt5}) and (\ref{Step2Sec4Eqt7}) that 
(\ref{Step2StatemSec4Eqt1}) is true. Step \ref{Step2ProofLemSec4} is proved.
\end{proof}

Step \ref{Step3ProofLemSec4} in the proof of (\ref{MainResToProve}) is concerned with the 
$L^1/L^{2^\star-1}$-controlled balance property of the $\tilde{\mathcal U}_\alpha$'s. 
Step \ref{Step3ProofLemSec4} is as follows. 

\begin{Step}\label{Step3ProofLemSec4} Let $\tilde{\mathcal U}_\alpha$ be 
given by (\ref{Sec4ProofEqt1}). There exists $C > 0$ such that, up to a subsequence,
 \begin{equation}\label{Step3StatemSec4Eqt1}
 \int_M\vert\tilde{\mathcal U}_\alpha\vert dv_g \le C 
 \int_M\vert\tilde{\mathcal U}_\alpha\vert^{2^\star-1}dv_g  
 \end{equation}
for all $\alpha$, where $\vert\tilde{\mathcal U}_\alpha\vert = \sum_{i=1}^p\vert\tilde u_\alpha^i\vert$, 
and $\vert\tilde{\mathcal U}_\alpha\vert^{2^\star-1} = \sum_{i=1}^p\vert\tilde u_\alpha^i\vert^{2^\star-1}$.
\end{Step}

\begin{proof}[Proof of Step \ref{Step3ProofLemSec4}] Let $f_\alpha^i = \hbox{sign}(\tilde u_\alpha^i)$ be the function given by
\begin{equation}\label{ProofStep3Sec4Eqt1}
f_\alpha^i = \chi_{\left\{\tilde u_\alpha^i > 0\right\}} - \chi_{\left\{\tilde u_\alpha^i < 0\right\}}\hskip.1cm ,
\end{equation}
where $\chi_A$ is the characterictic function of $A$. Then $f_\alpha^i\tilde u_\alpha^i = \vert\tilde u_\alpha^i\vert$ 
for all $\alpha$ and all $i$. We also have that $\vert f_\alpha^i\vert \le 1$ for all $\alpha$ and all $i$. 
As already mentionned, up to passing to a subsequence, 
we can assume that there exists $K > 0$ such that 
$A^\alpha_{ij}(x) \ge K\delta_{ij}$ in the sense of bilinear forms for all $x$. In particular, if we let $\Delta_g^p$ be the 
Laplacian acting on $p$-maps, the operators $\Delta_g^p + A(\alpha)$ are (uniformly) coercive in the sense that 
there exists $C > 0$ such that for any ${\mathcal U} \in H_{1,p}^2(M)$, and any $\alpha$,
\begin{equation}\label{ProofStep3Sec4Eqt2}
I_{A(\alpha)}({\mathcal U}) \ge C \Vert{\mathcal U}\Vert_{H_{1,p}^2}^2\hskip.1cm ,
\end{equation}
where
\begin{equation}\label{ProofStep3Sec4Eqt3}
I_{A(\alpha)}({\mathcal U}) = \int_M\vert\nabla{\mathcal U}\vert^2dv_g 
+ \int_MA(\alpha)({\mathcal U},{\mathcal U})dv_g\hskip.1cm .
\end{equation}
By (\ref{ProofStep3Sec4Eqt2}), and standard minimization technics, 
there is a solution ${\mathcal U}_\alpha^\prime$ to the minimization problem consisting of finding 
a minimizer for $I_{A(\alpha)}({\mathcal U})$ under the constraint 
$\int_M(f_\alpha,{\mathcal U})dv_g = 1$, where $I_{A(\alpha)}({\mathcal U})$ is as in 
(\ref{ProofStep3Sec4Eqt3}), 
$(f_\alpha,{\mathcal U}) = \sum_{i=1}^pf_\alpha^iu_i$, and the $u_i$'s are the components 
of ${\mathcal U}$. If $\lambda_\alpha$ is the minimum 
of the $I_{A(\alpha)}({\mathcal U})$'s, where ${\mathcal U} \in H_{1,p}^2(M)$ satisfies the constraint 
$\int_M(f_\alpha,{\mathcal U})dv_g = 1$, 
it easily follows from (\ref{ProofStep3Sec4Eqt2}) that $\lambda_\alpha > 0$. 
We let $\hat{\mathcal U}_\alpha = \lambda_\alpha^{-1}{\mathcal U}_\alpha^\prime$. 
 Then $\hat{\mathcal U}_\alpha$ is a solution of the system
\begin{equation}\label{ProofStep3Sec4Eqt4}
\Delta_g\hat u_\alpha^i + \sum_{j=1}^pA^\alpha_{ij}\hat u_\alpha^j = f_\alpha^i
\end{equation}
for all $i$ and all $\alpha$, 
where the $\hat u_\alpha^i$'s are the components of $\hat{\mathcal U}_\alpha$, and the 
$f_\alpha^i$'s are as in (\ref{ProofStep3Sec4Eqt1}). 
Multiplying (\ref{ProofStep3Sec4Eqt4}) by $\hat u_\alpha^i$, integrating over $M$, and summing over $i$, 
we get with (\ref{ProofStep3Sec4Eqt2}) that the square of the 
$H_{1,p}^2$-norm of the $\hat{\mathcal U}_\alpha$'s is uniformly controlled 
by the $L^1$-norm of the $\vert\hat{\mathcal U}_\alpha\vert$'s. 
In particular, the $\hat u_\alpha^i$'s are uniformly bounded in $L^2$. By standard elliptic 
theory, the $\hat u_\alpha^i$'s are in the 
Sobolev spaces $H_2^q$ for all $q$. As an easy consequence, 
the $\hat u_\alpha^i$'s are continuous. By the above discussion, and standard elliptic theory, we then  
get that there exists $C > 0$ such that $\vert\hat u_\alpha^i\vert \le C$ in $M$, for all $\alpha$ and all $i$. 
By (\ref{Sec4ProofEqt2}) and (\ref{ProofStep3Sec4Eqt4}) we can now write that
\begin{equation}\label{ProofStep3Sec4Eqt5}
\begin{split}
\sum_{i=1}^p\int_M\vert\tilde u_\alpha^i\vert dv_g
&= \sum_{i=1}^p\int_M\tilde u_\alpha^if_\alpha^idv_g\\
&= \sum_{i=1}^p\int_M\left(\Delta_g\hat u_\alpha^i + \sum_{j=1}^pA^\alpha_{ij}\hat u_\alpha^j\right)\tilde u_\alpha^idv_g\\
&= \sum_{i=1}^p\int_M\left(\Delta_g\tilde u_\alpha^i 
+ \sum_{j=1}^pA^\alpha_{ij}\tilde u_\alpha^j\right)\hat u_\alpha^idv_g\\
&= \sum_{i=1}^p\int_M\vert\tilde u_\alpha^i\vert^{2^\star-2}\tilde u_\alpha^i\hat u_\alpha^idv_g\\
\end{split}
\end{equation}
for all $\alpha$. Since there exists $C > 0$ such that $\vert\hat u_\alpha^i\vert \le C$ in $M$ 
for all $\alpha$ and all $i$, it follows from (\ref{ProofStep3Sec4Eqt5}) that
$$\int_M\vert\tilde{\mathcal U}_\alpha\vert dv_g \le C \int_M\vert\tilde{\mathcal U}_\alpha\vert^{2^\star-1}dv_g$$
for all $\alpha$, 
where $C > 0$ does not depend on $\alpha$. This proves Step \ref{Step3ProofLemSec4}.
\end{proof}

Step \ref{Step4ProofLemSec4} in the proof of (\ref{MainResToProve}) is concerned with 
$L^2$-concentration. 
We assume here that $n \ge 4$. When $n =3$, bubbles do not concentrate in the $L^2$-norm,  
and $L^2$-concentration turns out to be false in this dimension. Dimension $4$ is the 
smallest dimension for this notion of $L^2$-concentration.
Step \ref{Step4ProofLemSec4} is as follows. 

\begin{Step}\label{Step4ProofLemSec4} Let $\tilde{\mathcal U}_\alpha$ and $\tilde{\mathcal U}^0$ be 
given by (\ref{Sec4ProofEqt1}) and (\ref{Sec4ProofEqt4}). Assume $\tilde{\mathcal U}^0 \equiv 0$ and $n \ge 4$. 
 Up to a subsequence,
 \begin{equation}\label{Step4StatemSec4Eqt1}
 \lim_{\alpha \to +\infty} 
 \frac{\int_{B_\delta}\vert\tilde{\mathcal U}_\alpha\vert^2dv_g}
 {\int_M\vert\tilde{\mathcal U}_\alpha\vert^2dv_g} 
 = 1
 \end{equation}
for all $\delta > 0$, where $B_\delta = B_{x_0}(\delta)$ is the ball centered at $x_0$ of 
radius $\delta$, 
$x_0$ is the limit of the centers of the $1$-bubble from which 
the $p$-bubble $({\mathcal B}_\alpha)_\alpha$ in (\ref{Step1StatemSec4Eqt1}) is defined, 
and $\vert\tilde{\mathcal U}_\alpha\vert^2 = \sum_{i=1}^p\vert\tilde u_\alpha^i\vert^2$.
\end{Step}

\begin{proof}[Proof of Step \ref{Step4ProofLemSec4}] Clearly, Step \ref{Step4ProofLemSec4} is equivalent 
to proving that for any $\delta > 0$, $R_\delta(\alpha) \to 0$ as $\alpha \to +\infty$, where $R_\delta(\alpha)$ is the ratio 
given by
\begin{equation}\label{RatioL2Conc}
R_\delta(\alpha) = \frac{\int_{M\backslash B_\delta}\vert{\mathcal U}_\alpha\vert^2dv_g}
{\int_M\vert{\mathcal U}_\alpha\vert^2dv_g}\hskip.1cm .
\end{equation}
We fix $\delta > 0$. By Steps \ref{Step2ProofLemSec4} and \ref{Step3ProofLemSec4}, 
we can write that for any $\alpha$, 
\begin{eqnarray*}
\int_{M\backslash{\mathcal B}_\delta}\vert\tilde{\mathcal U}_\alpha\vert^2dv_g
& \le & \left(\max_{M\backslash{\mathcal B}_\delta}\vert\tilde{\mathcal U}_\alpha\vert\right) 
\int_{M\backslash{\mathcal B}_\delta}\vert\tilde{\mathcal U}_\alpha\vert dv_g\\
& \le & C \sqrt{\int_M\vert\tilde{\mathcal U}_\alpha\vert^2dv_g} 
\int_M\vert\tilde{\mathcal U}_\alpha\vert^{2^\star-1}dv_g\hskip.1cm ,
\end{eqnarray*}
where $C > 0$ is independent of $\alpha$. In particular, 
\begin{equation}\label{RatioL2ConcProof}
R_\delta(\alpha) \le C \frac{\int_M\vert\tilde{\mathcal U}_\alpha\vert^{2^\star-1}dv_g}
{\sqrt{\int_M\vert\tilde{\mathcal U}_\alpha\vert^2dv_g}}
\end{equation}
for all $\alpha$, where $C > 0$ is independent of $\alpha$, 
and $R_\delta(\alpha)$ is given by (\ref{RatioL2Conc}). If we assume now that $n \ge 6$, then 
$2^\star-1 \le 2$, and we can write with H\"older's inequality that
$$\int_M\vert\tilde u_\alpha^i\vert^{2^\star-1}dv_g \le V_g^{\frac{3-2^\star}{2}} 
\left(\int_M\vert\tilde u_\alpha^i\vert^2dv_g\right)^{\frac{2^\star-1}{2}}$$
for all $i$, where $V_g$ is the volume of $M$ with respect to $g$. 
In particular, there exists 
$C > 0$ such that
\begin{equation}\label{HighDimProofStep4}
\int_M\vert\tilde{\mathcal U}_\alpha\vert^{2^\star-1}dv_g \le C 
\left(\int_M\vert\tilde{\mathcal U}_\alpha\vert^2dv_g\right)^{\frac{2^\star-1}{2}}\hskip.1cm ,
\end{equation}
for all $\alpha$. 
Since $\tilde{\mathcal U}^0 \equiv 0$, it follows from (\ref{Sec4ProofEqt4}) that 
$\tilde{\mathcal U}_\alpha \to 0$ in $L^2$ as $\alpha \to +\infty$. 
Since $2^\star > 2$, we then get with 
(\ref{RatioL2ConcProof}) and (\ref{HighDimProofStep4}) that $R_\delta(\alpha) \to 0$ as 
$\alpha \to +\infty$. 
This proves (\ref{Step4StatemSec4Eqt1}) when $n \ge 6$. 
If we assume now that $n = 5$, then 
$2 \le 2^\star-1 \le 2^\star$, and we can write with H\"older's inequality that
$$\left(\int_M\vert\tilde u_\alpha^i\vert^{2^\star-1}dv_g\right)^{\frac{1}{2^\star-1}}
\le \left(\int_M\vert\tilde u_\alpha^i\vert^2dv_g\right)^{\frac{\theta}{2}}
\left(\int_M\vert\tilde u_\alpha^i\vert^{2^\star}dv_g\right)^{\frac{1-\theta}{2^\star}}\hskip.1cm ,$$
where $\theta = \frac{3}{2(2^\star-1)}$. By (\ref{Sec4ProofEqt3}) we then get that
$$\int_M\vert\tilde{\mathcal U}_\alpha\vert^{2^\star-1}dv_g \le C 
\left(\int_M\vert\tilde{\mathcal U}_\alpha\vert^2dv_g\right)^{\frac{3}{4}}\hskip.1cm ,$$
where $C > 0$ does not depend on $\alpha$. 
Since $\frac{3}{4} > \frac{1}{2}$ and 
$\tilde{\mathcal U}_\alpha \to 0$ in $L^2$, we get with 
(\ref{RatioL2ConcProof}) that $R_\delta(\alpha) \to 0$ as 
$\alpha \to +\infty$. This proves (\ref{Step4StatemSec4Eqt1}) when $n = 5$. 
Now it remains to prove 
(\ref{Step4StatemSec4Eqt1}) when $n = 4$. The argument when $n = 4$ is slightly more delicate. 
We start writing that
\begin{equation}\label{L2ConcProofDim4Eqt1}
\begin{split}
\frac{\int_M\vert\tilde{\mathcal U}_\alpha\vert^{2^\star-1}dv_g}{\sqrt{\int_M\vert\tilde{\mathcal U}_\alpha\vert^2dv_g}} 
&= \sum_{i=1}^p
\frac{\int_M\vert\tilde u_\alpha^i\vert^{2^\star-1}dv_g}{\sqrt{\int_M\vert\tilde{\mathcal U}_\alpha\vert^2dv_g}}\\
&\le 
\sum_{i=1}^p\frac{\int_M\vert\tilde u_\alpha^i\vert^{2^\star-1}dv_g}{\sqrt{\int_M\vert\tilde u_\alpha^i\vert^2dv_g}}\hskip.1cm .
\end{split}
\end{equation}
We let  the $x_\alpha$'s and $\mu_\alpha$'s be the centers and weights 
of the $1$-bubble $(B_\alpha)_\alpha$ 
from which the $p$-bubble $({\mathcal B}_\alpha)_\alpha$ in (\ref{Step1StatemSec4Eqt1}) is defined. 
We let $i_0 = 1,\dots,p$, be such that ${\mathcal B}_\alpha^{i_0} = B_\alpha$ for all $\alpha$. For $R > 0$, 
we let also $\Omega_{i_0,\alpha}(R)$ be given by 
\begin{equation}\label{L2ConcProofDim4EqtSet1}
\Omega_{i_0,\alpha}(R) = B_{x_\alpha}(R\mu_\alpha)\hskip.1cm .
\end{equation}
Since $n = 4$, we have that $2^\star = 4$. If $i \not= i_0$, we can write, thanks to H\"older's inequalities, that
$$\int_M\vert\tilde u_\alpha^i\vert^{2^\star-1}dv_g 
\le \sqrt{\int_M\vert\tilde u_\alpha^i\vert^{2^\star}dv_g}
\sqrt{\int_M\vert\tilde u_\alpha^i\vert^2dv_g}\hskip.1cm .$$
By (\ref{Step1StatemSec4Eqt1}), when $i \not= i_0$, $\tilde u_\alpha^i \to 0$ in $H_1^2(M)$ as $\alpha \to +\infty$. 
It follows that for any $i \not= i_0$, 
\begin{equation}\label{Sec4ProofL2ConcConclEqt1}
\frac{\int_M\vert\tilde u_\alpha^i\vert^{2^\star-1}dv_g}{\sqrt{\int_M\vert\tilde u_\alpha^i\vert^2dv_g}} = o(1)\hskip.1cm ,
\end{equation}
where $o(1) \to 0$ as $\alpha \to +\infty$.
On the other hand, when $i = i_0$, we get with H\"older's inequalities that
\begin{eqnarray*}
\int_M\vert\tilde u_\alpha^{i_0}\vert^{2^\star-1}dv_g 
&\le& \int_{\Omega_{i_0,\alpha}(R)}\vert\tilde u_\alpha^{i_0}\vert^{2^\star-1}dv_g\\
&&\hskip.4cm + \sqrt{\int_{M\backslash\Omega_{i_0,\alpha}(R)}\vert\tilde u_\alpha^{i_0}\vert^{2^\star}dv_g}
\sqrt{\int_M\vert\tilde u_\alpha^{i_0}\vert^2dv_g}\hskip.1cm ,
\end{eqnarray*}
and we can write that
\begin{equation}\label{L2ConcProofDim4Eqt3}
\frac{\int_M\vert\tilde u_\alpha^{i_0}\vert^{2^\star-1}dv_g}{\sqrt{\int_M\vert\tilde u_\alpha^{i_0}\vert^2dv_g}} 
\le \sqrt{\int_{M\backslash\Omega_{i_0,\alpha}(R)}\vert\tilde u_\alpha^{i_0}\vert^{2^\star}dv_g} 
+ \frac{\int_{\Omega_{i_0,\alpha}(R)}\vert\tilde u_\alpha^{i_0}\vert^{2^\star-1}dv_g}
{\sqrt{\int_M\vert\tilde u_\alpha^{i_0}\vert^2dv_g}}\hskip.1cm ,
\end{equation}
where $\Omega_{i_0,\alpha}(R)$ is as in (\ref{L2ConcProofDim4EqtSet1}). 
For $\varphi \in C^\infty_0({\mathbb R}^n)$, where $C^\infty_0({\mathbb R}^n)$ is the set of smooth functions with 
compact support in ${\mathbb R}^n$, we let $\varphi_\alpha^{i_0}$ be
the function defined by the equation 
\begin{equation}\label{L2ConcProofDim4Eqt4}
\varphi_\alpha^{i_0}(x) = (\mu_\alpha)^{-\frac{n-2}{2}}
\varphi\left((\mu_\alpha)^{-1}\exp_{x_\alpha}^{-1}(x)\right)
\hskip.1cm .
\end{equation}
Straightforward computations give that 
for any $R > 0$,\par
\medskip (i) $\displaystyle\int_{M\backslash \Omega_{i_0,\alpha}(R)}
(B_\alpha)^{2^\star}dv_g = \varepsilon_R(\alpha)$,\par
\medskip (ii) $\displaystyle\int_{\Omega_{i_0,\alpha}(R)}
(B_\alpha)^{2^\star-1}\varphi_\alpha^{i_0}dv_g 
= \int_{B_0(R)}u^{2^\star-1}\varphi dx + o(1)$,\par
\medskip (iii) $\displaystyle\int_{\Omega_{i_0,\alpha}(R)}
(B_\alpha)^2(\varphi_\alpha^{i_0})^{2^\star-2}dv_g 
= \int_{B_0(R)}u^2\varphi^{2^\star-2}dx + o(1)$\par
\medskip\noindent where $i_0$ is such that ${\mathcal B}_\alpha^{i_0} = B_\alpha$ for all $\alpha$, 
$({\mathcal B}_\alpha)_\alpha$ is the $p$-bubble in (\ref{Step1StatemSec4Eqt1}), 
$u$ is given by (\ref{PosSolCritEuclEqt}), $\Omega_{i_0,\alpha}(R)$ is given by (\ref{L2ConcProofDim4EqtSet1}), 
$o(1) \to 0$ as $\alpha \to +\infty$, and the $\varepsilon_R(\alpha)$'s are such that
\begin{equation}\label{EqtRestL2ConcDim4}
\lim_{R \to +\infty}\limsup_{\alpha\to +\infty}\varepsilon_R(\alpha) = 0
\hskip.1cm .
\end{equation}
By (i) and (\ref{Step1StatemSec4Eqt1}) we can write that
\begin{equation}\label{L2ConcProofDim4Eqt5}
\int_{M\backslash\Omega_{i_0,\alpha}(R)}\vert\tilde u_\alpha^{i_0}\vert^{2^\star}dv_g = \varepsilon_R(\alpha)\hskip.1cm ,
\end{equation}
where $\Omega_{i_0,\alpha}(R)$ is as in (\ref{L2ConcProofDim4EqtSet1}), and 
the $\varepsilon_R(\alpha)$'s are such that (\ref{EqtRestL2ConcDim4}) holds.
From now on, we let $\varphi$ in (\ref{L2ConcProofDim4Eqt4}) be such that $\varphi = 1$ in 
$B_0(R)$, $R > 0$. Then,
$$\int_{\Omega_{i_0,\alpha}(R)}\vert\tilde u_\alpha^{i_0}\vert^{2^\star-1}dv_g 
= \mu_\alpha^{\frac{n-2}{2}}
\int_{\Omega_{i_0,\alpha}(R)}\vert\tilde u_\alpha^{i_0}\vert^{2^\star-1}\varphi_\alpha^{i_0}dv_g$$
and, by  (\ref{Step1StatemSec4Eqt1}) and (ii), we can write that
\begin{eqnarray*} \int_{\Omega_{i_0,\alpha}(R)}\vert\tilde u_\alpha^{i_0}\vert^{2^\star-1}\varphi_\alpha^{i_0}dv_g
& \le & C\int_{\Omega_{i_0,\alpha}(R)}B_\alpha^{2^\star-1}\varphi_\alpha^{i_0}dv_g + o(1)\\
& \le & C \int_{B_0(R)}u^{2^\star-1}dx + o(1)\hskip.1cm ,
\end{eqnarray*}
where $o(1) \to 0$ as $\alpha \to +\infty$, and $C > 0$ does not depend on $\alpha$ 
and $R$. In particular, we have that
\begin{equation}\label{L2ConcProofDim4Eqt6}
\int_{\Omega_{i_0,\alpha}(R)}\vert\tilde u_\alpha^{i_0}\vert^{2^\star-1}dv_g 
\le \left(C \int_{B_0(R)}u^{2^\star-1}dx + o(1)\right)
\mu_\alpha^{\frac{n-2}{2}}\hskip.1cm ,
\end{equation}
where $o(1) \to 0$ as $\alpha \to +\infty$, $u$ is as in (\ref{PosSolCritEuclEqt}), 
and $C > 0$ does not depend on $\alpha$ and $R$.
Independently, we also have that
\begin{eqnarray*} \int_M\vert\tilde u_\alpha^{i_0}\vert^2dv_g
& \ge & \int_{\Omega_{i_0,\alpha}(R)}\vert\tilde u_\alpha^{i_0}\vert^2dv_g\\
& \ge & \mu_\alpha^{n-2}
\int_{\Omega_{i_0,\alpha}(R)}\vert\tilde u_\alpha^{i_0}\vert^2(\varphi_\alpha^{i_0})^{2^\star-2}dv_g
\end{eqnarray*}
Here, $2^\star-2 = 2$. As is easily checked, we can write 
with (\ref{Step1StatemSec4Eqt1}) that
$$\int_{\Omega_{i_0,\alpha}(R)}\vert\tilde u_\alpha^{i_0}\vert^2(\varphi_\alpha^{i_0})^{2^\star-2}dv_g
 =  \int_{\Omega_{i_0,\alpha}(R)}B_\alpha^2
(\varphi_\alpha^{i_0})^{2^\star-2}dv_g + o(1)$$
and thanks to (iii) we get that
$$\int_{\Omega_{i_0,\alpha}(R)}\vert\tilde u_\alpha^{i_0}\vert^2(\varphi_\alpha^{i_0})^{2^\star-2}dv_g 
\ge \int_{B_0(R)}u^2dx + o(1)\hskip.1cm .$$
In particular,
\begin{equation}\label{L2ConcProofDim4Eqt7}
\int_M\vert\tilde u_\alpha^{i_0}\vert^2dv_g 
\ge \mu_\alpha^{n-2} \left(\int_{B_0(R)}u^2dx + o(1)\right)\hskip.1cm ,
\end{equation}
where $o(1) \to 0$ as $\alpha \to +\infty$, 
and $u$ is as in (\ref{PosSolCritEuclEqt}). By (\ref{RatioL2ConcProof}), (\ref{L2ConcProofDim4Eqt1}), and (\ref{Sec4ProofL2ConcConclEqt1}),
we can write that
$$R_\delta(\alpha) 
\le C \frac{\int_M\vert\tilde u_\alpha^{i_0}\vert^{2^\star-1}dv_g}{\sqrt{\int_M\vert\tilde u_\alpha^{i_0}\vert^2dv_g}} 
+ o(1)$$
for all $\alpha$, where $R_\delta(\alpha) $ is given 
by (\ref{RatioL2Conc}), and $C > 0$ is independent of $\alpha$. Then, by 
(\ref{L2ConcProofDim4Eqt3}), and (\ref{L2ConcProofDim4Eqt5})--(\ref{L2ConcProofDim4Eqt7}), we get that
for any $R > 0$, 
\begin{equation}\label{L2ConcProofDim4Eqt8}
\limsup_{\alpha\to+\infty}R_\delta(\alpha) \le \varepsilon_R + 
C \frac{\int_{B_0(R)}u^{2^\star-1}dx}{\sqrt{\int_{B_0(R)}u^2dx}}\hskip.1cm ,
\end{equation}
where $\varepsilon_R \to 0$ as $R \to +\infty$, and $C > 0$ does not depend 
on $R$. It is easily seen that
\begin{eqnarray*} \lim_{R\to+\infty}\int_{B_0(R)}u^{2^\star-1}dx
& = & \int_{{\mathbb R}^n}u^{2^\star-1}dx\\
& < & +\infty
\end{eqnarray*}
On the other hand, when $n = 4$,
$$\lim_{R\to+\infty}\int_{B_0(R)}u^2dx = +\infty\hskip.1cm .$$
Coming back to (\ref{L2ConcProofDim4Eqt8}), it follows that for any $\delta > 0$, 
$R_\delta(\alpha) \to 0$ as $\alpha \to +\infty$. In particular, (\ref{Step4StatemSec4Eqt1}) 
is true when $n = 4$.  This ends the proof of Step \ref{Step4ProofLemSec4}.
\end{proof}

Step \ref{Step5ProofLemSec4} in the proof of (\ref{MainResToProve}) is concerned with 
proving that the off diagonal terms $\int_M\tilde u_\alpha^i\tilde u_\alpha^jdv_g$, $i \not= j$, are small 
when compared to the diagonal terms $\int_M(\tilde u_\alpha^i)^2dv_g$. 
Step \ref{Step5ProofLemSec4} is as follows.

\begin{Step}\label{Step5ProofLemSec4} Let $\tilde{\mathcal U}_\alpha$ and $\tilde{\mathcal U}^0$ be 
given by (\ref{Sec4ProofEqt1}) and (\ref{Sec4ProofEqt4}). Assume $\tilde{\mathcal U}^0 \equiv 0$. 
Up to a subsequence, for any $i, j = 1,\dots,p$ such that $i\not= j$,  
 \begin{equation}\label{Step5StatemSec4Eqt1}
\int_{B_{x_0}(\delta)}\vert\tilde u_\alpha^i\tilde u_\alpha^j\vert dv_g 
\le \varepsilon_\delta \int_M\vert\tilde{\mathcal U}_\alpha\vert^2dv_g
 \end{equation}
for all $\delta > 0$ and all $\alpha$, 
where $x_0$ is the limit of the centers of the $1$-bubble from which 
the $p$-bubble $({\mathcal B}_\alpha)_\alpha$ in (\ref{Step1StatemSec4Eqt1}) is defined, 
$\vert\tilde{\mathcal U}_\alpha\vert^2 = \sum_{i=1}^p\vert\tilde u_\alpha^i\vert^2$, and 
$\varepsilon_\delta > 0$ is independent of $\alpha$ and such that 
$\varepsilon_\delta \to 0$ as $\delta \to 0$. 
\end{Step}

\begin{proof}[Proof of Step \ref{Step5ProofLemSec4}] As in the proof 
of Step \ref{Step4ProofLemSec4}, we let $i_0 = 1,\dots,p$, be such that 
${\mathcal B}_\alpha^{i_0} = B_\alpha$ for all $\alpha$, 
where $(B_\alpha)_\alpha$ is the $1$-bubble 
from which the $p$-bubble $({\mathcal B}_\alpha)_\alpha$ in (\ref{Step1StatemSec4Eqt1}) is defined.
Then, by (\ref{Step1StatemSec4Eqt1}),
\begin{equation}\label{ProofStep5Sec4Eqt1}
\int_M\vert\tilde u_\alpha^i\vert^{2^\star}dv_g = o(1)
\end{equation}
for all $\alpha$ and all $i \not= i_0$, where $o(1) \to 0$ as $\alpha \to +\infty$. Let $i\not= i_0$. 
We multiply the $i$th equation in (\ref{Sec4ProofEqt2}) by $\tilde u_\alpha^i$, and integrate over $M$. Then
 we can write that
 \begin{equation}\label{ProofStep5Sec4Eqt2}
 \int_M\left(\vert\nabla\tilde u_\alpha^i\vert^2 + A^\alpha_{ii}(\tilde u_\alpha^i)^2\right)dv_g 
 \le \int_M\vert\tilde u_\alpha^i\vert^{2^\star}dv_g + 
 C \sum_{j\not= i}\int_M\vert\tilde u_\alpha^i\vert \vert\tilde u_\alpha^j\vert dv_g
 \end{equation}
 for all $\alpha$, where $C > 0$ is independent of $\alpha$ and $i$. As already mentionned in the 
 introduction of this section, 
 up to passing to a subsequence, we can assume that 
 there exists $K > 0$ such that $A^\alpha_{ij} \ge K\delta_{ij}$ in the sense of bilinear 
 forms, for all $\alpha$. Then $A^\alpha_{ii} \ge K$ in $M$, for all $\alpha$ and all $i$, and by the 
 Sobolev embedding theorem, we get that there exists $C > 0$ such that
 \begin{equation}\label{ProofStep5Sec4Eqt3}
 \int_M\left(\vert\nabla\tilde u_\alpha^i\vert^2 + A^\alpha_{ii}(\tilde u_\alpha^i)^2\right)dv_g 
 \ge C \left(\int_M\vert\tilde u_\alpha^i\vert^{2^\star}dv_g\right)^{2/2^\star}
 \end{equation}
 for all $\alpha$. Combining (\ref{ProofStep5Sec4Eqt2}) and (\ref{ProofStep5Sec4Eqt3}), we get that 
 there exist $C, C^\prime > 0$ such that
  \begin{equation}\label{ProofStep5Sec4Eqt4}
 C \left(\int_M\vert\tilde u_\alpha^i\vert^{2^\star}dv_g\right)^{2/2^\star}
 \le \int_M\vert\tilde u_\alpha^i\vert^{2^\star}dv_g + 
 C^\prime \sum_{j\not= i}\int_M\vert\tilde u_\alpha^i\vert \vert\tilde u_\alpha^j\vert dv_g
 \end{equation}
 for all $\alpha$, and all $i\not= i_0$. 
 By H\"older's inequality,
 \begin{equation}\label{ProofStep5Sec4Eqt5}
 \begin{split} \int_M\vert\tilde u_\alpha^i\vert \vert\tilde u_\alpha^j\vert dv_g
 &\le \sqrt{\int_M\vert\tilde u_\alpha^i\vert^2dv_g} \sqrt{\int_M\vert\tilde{\mathcal U}_\alpha\vert^2dv_g}\\
 &\le C\left(\int_M\vert\tilde u_\alpha^i\vert^{2^\star}dv_g\right)^{1/2^\star} 
 \sqrt{\int_M\vert\tilde{\mathcal U}_\alpha\vert^2dv_g}
\end{split}
\end{equation}
for all $\alpha$, where $C > 0$ is independent of $\alpha$. Combining 
(\ref{ProofStep5Sec4Eqt4}) and (\ref{ProofStep5Sec4Eqt5}), it follows that
 \begin{equation}\label{ProofStep5Sec4Eqt6}
 \begin{split}
 C \left(\int_M\vert\tilde u_\alpha^i\vert^{2^\star}dv_g\right)^{2/2^\star}
& \le \int_M\vert\tilde u_\alpha^i\vert^{2^\star}dv_g\\
&\hskip.4cm 
+ C^\prime \left(\int_M\vert\tilde u_\alpha^i\vert^{2^\star}dv_g\right)^{1/2^\star}
 \sqrt{\int_M\vert\tilde{\mathcal U}_\alpha\vert^2dv_g}
 \end{split}
 \end{equation}
for all $\alpha$, and all $i\not= i_0$, where $C, C^\prime > 0$ are independent of $\alpha$ and $i$. By 
(\ref{ProofStep5Sec4Eqt1}) we then get that there exists $C > 0$ such that
 \begin{equation}\label{ProofStep5Sec4Eqt7}
\left(\int_M\vert\tilde u_\alpha^i\vert^{2^\star}dv_g\right)^{1/2^\star} \le C 
\sqrt{\int_M\vert\tilde{\mathcal U}_\alpha\vert^2dv_g}
\end{equation}
for all $\alpha$, and all $i\not= i_0$. Now, given $\delta > 0$, and $i\not= j$ arbitrary, we write that
\begin{equation}\label{ProofStep5Sec4Eqt8}
\int_{B_{x_0}(\delta)}\vert\tilde u_\alpha^i\vert \vert\tilde u_\alpha^j\vert dv_g
\le \sqrt{\int_{B_{x_0}(\delta)}(\tilde u_\alpha^i)^2dv_g}
\sqrt{\int_{B_{x_0}(\delta)}(\tilde u_\alpha^j)^2dv_g}
\end{equation}
for all $\alpha$, where $x_0$ is the limit of the centers of the $1$-bubble from which 
the $p$-bubble $({\mathcal B}_\alpha)_\alpha$ in (\ref{Step1StatemSec4Eqt1}) is defined. 
Since $i\not= j$, either $i\not= i_0$ or $j\not= i_0$. Suppose $j \not= i_0$. On the one 
hand we can write that
\begin{equation}\label{ProofStep5Sec4Eqt9}
\int_{B_{x_0}(\delta)}(\tilde u_\alpha^i)^2dv_g\
\le \int_M\vert\tilde{\mathcal U}_\alpha\vert^2dv_g
\end{equation}
for all $\alpha$ and $\delta > 0$. 
On the other hand, by H\"older's inequality, we can write that
\begin{eqnarray*}
\int_{B_{x_0}(\delta)}(\tilde u_\alpha^j)^2dv_g
&&\le \vert B_{x_0}(\delta)\vert^{\frac{2^\star-2}{2^\star}}
\left(\int_{B_{x_0}(\delta)}\vert\tilde u_\alpha^j\vert^{2^\star}dv_g\right)^{2/2^\star}\\
&&\le \vert B_{x_0}(\delta)\vert^{\frac{2^\star-2}{2^\star}}
\left(\int_M\vert\tilde u_\alpha^j\vert^{2^\star}dv_g\right)^{2/2^\star}
\end{eqnarray*}
for all $\alpha$, where $\vert B_{x_0}(\delta)\vert$ is the volume of $B_{x_0}(\delta)$ 
with respect to $g$. By  (\ref{ProofStep5Sec4Eqt7}), since $j \not= i_0$, we then get that
\begin{equation}\label{ProofStep5Sec4Eqt10}
\int_{B_{x_0}(\delta)}(\tilde u_\alpha^j)^2dv_g \le 
C \vert B_{x_0}(\delta)\vert^{\frac{2^\star-2}{2^\star}} \int_M\vert\tilde{\mathcal U}_\alpha\vert^2dv_g
\end{equation}
for all $\alpha$ and $\delta > 0$, where $C > 0$ is independent of $\alpha$ and $\delta$. 
Plugging (\ref{ProofStep5Sec4Eqt9}) and (\ref{ProofStep5Sec4Eqt10}) into (\ref{ProofStep5Sec4Eqt8}), 
since $\vert B_{x_0}(\delta)\vert \to 0$ as $\delta \to 0$, 
we get that (\ref{Step5StatemSec4Eqt1}) is true. This ends the proof of Step \ref{Step5ProofLemSec4}.
\end{proof}

By Steps \ref{Step1ProofLemSec4} to \ref{Step5ProofLemSec4}, we are now in position to prove 
(\ref{MainResToProve}), and hence to prove Lemma \ref{TheLem}. We use in the process that 
for any $\varepsilon > 0$, there 
exists $\delta_\varepsilon > 0$, such that for any smooth function $u$ with compact support 
in $B_{x_0}(\delta_\varepsilon)$,
\begin{equation}\label{BOLoc}
\left(\int_M\vert u\vert^{2^\star}dv_g\right)^{2/2^\star} 
\le K_n^2\int_M\vert\nabla u\vert_{\hat g}^2dv_g + B_\varepsilon\int_Mu^2dv_g
\end{equation}
where $B_\varepsilon = \frac{n-2}{4(n-1)}K_n^2\left(S_g(x_0) + \varepsilon\right)$, 
$K_n$ is given by (\ref{SharpCstEucl}), and 
$S_g$ is the scalar curvature of $g$. 
Inequality (\ref{BOLoc}) is a straightforward consequence of the 
local isoperimetric inequality proved in Druet \cite{Dru}. Step \ref{Step6ProofLemSec4} is as follows.

\begin{Step}\label{Step6ProofLemSec4} Let $\tilde{\mathcal U}_\alpha$ and $\tilde{\mathcal U}^0$ be given by 
(\ref{Sec4ProofEqt1}) and (\ref{Sec4ProofEqt4}). Assume that for any $i$ and any $x \in M$,
\begin{equation}\label{IneqSecProof4Eqt1}
A^0_{ii}(x) > \frac{n-2}{4(n-1)}S_g(x)\hskip.1cm ,
\end{equation}
and that $n \ge 4$. Then 
$\tilde{\mathcal U}^0 \not\equiv 0$. In particular, (\ref{MainResToProve}) and Lemma \ref{TheLem} are true.
\end{Step}

\begin{proof}[Proof of Step \ref{Step6ProofLemSec4}] We proceed by contradiction 
and assume that $\tilde{\mathcal U}^0 \equiv 0$. We let 
$x_0$ be the limit of the centers of the $1$-bubble from which 
the $p$-bubble $({\mathcal B}_\alpha)_\alpha$ in (\ref{Step1StatemSec4Eqt1}) is defined. 
We fix $\varepsilon > 0$, 
and let $\eta$ be a smooth cutoff function such that $\eta = 1$ in 
$B_{x_0}(\delta_\varepsilon/4)$, $\eta = 0$ in 
$M\backslash B_{x_0}(\delta_\varepsilon/2)$, and $0 \le \eta \le 1$. We
plugg the $\eta\tilde u_\alpha^i$'s into (\ref{BOLoc}), $i = 1,\dots,p$, and then sum over $i$. 
Noting that
$$\int_M\vert\nabla(\eta\tilde u_\alpha^i)\vert^2dv_g = 
\int_M\eta^2\tilde u_\alpha^i(\Delta_g\tilde u_\alpha^i)dv_g 
+ \int_M\vert\nabla\eta\vert^2(\tilde u_\alpha^i)^2dv_g\hskip.1cm ,$$
and that $\vert\nabla\eta\vert = 0$ around $x_0$, 
we get with (\ref{Sec4ProofEqt2}) and $L^2$-concentration in Step \ref{Step4ProofLemSec4} that
\begin{equation}\label{Step6Sec4ProofEqt1}
\begin{split}
&\sum_{i=1}^p\left(\left(\int_M\vert\eta\tilde u_\alpha^i\vert^{2^\star}dv_g\right)^{2/2^\star} 
- K_n^2\int_M\eta^2\vert\tilde u_\alpha^i\vert^{2^\star}dv_g\right)\\
&\le -K_n^2\sum_{i,j=1}^p \int_M\eta^2A_{ij}^\alpha\tilde u_\alpha^i\tilde u_\alpha^jdv_g 
+ \left(B_\varepsilon + o(1)\right)\int_M\vert\tilde{\mathcal U}_\alpha\vert^2dv_g
\end{split}
\end{equation}
for all $\alpha$, where $o(1) \to 0$ as $\alpha \to +\infty$, and $B_\varepsilon$ is as in (\ref{BOLoc}). 
By H\"older's inequality, and (\ref{Sec4ProofEqt3}),
\begin{equation}\label{Step6Sec4ProofEqt2}
\begin{split}
\int_M\eta^2\vert\tilde u_\alpha^i\vert^{2^\star}dv_g
&\le \left(\int_M\vert\eta\tilde u_\alpha^i\vert^{2^\star}dv_g\right)^{2/2^\star}
\left(\int_M\vert\tilde u_\alpha^i\vert^{2^\star}dv_g\right)^{(2^\star-2)/2^\star}\\
&\le K_n^{-2}\left(\int_M\vert\eta\tilde u_\alpha^i\vert^{2^\star}dv_g\right)^{2/2^\star}
\end{split}
\end{equation}
for all $\alpha$ and $i$. By (\ref{Step6Sec4ProofEqt2}), the left hand side in 
(\ref{Step6Sec4ProofEqt1}) is nonnegative. Since we also have that 
$A^\alpha_{ij} \to A^0_{ij}$ in $C^{0,\theta}(M)$, 
we can write with (\ref{Step6Sec4ProofEqt1}) that
\begin{equation}\label{Step6Sec4ProofEqt3}
K_n^2\sum_{i,j=1}^p \int_M\eta^2A_{ij}^0\tilde u_\alpha^i\tilde u_\alpha^jdv_g 
\le \left(B_\varepsilon + o(1)\right)\int_M\vert\tilde{\mathcal U}_\alpha\vert^2dv_g
\end{equation}
for all $\alpha$, where $o(1) \to 0$ as $\alpha \to +\infty$, and $B_\varepsilon$ is as in (\ref{BOLoc}). 
By $L^2$-concentration in Step \ref{Step4ProofLemSec4}, and the control of 
the off diagonal terms in Step \ref{Step5ProofLemSec4}, 
it follows from (\ref{Step6Sec4ProofEqt3}) that for any $\varepsilon > 0$, and any $\delta > 0$, 
\begin{equation}\label{Step6Sec4ProofEqt4}
\begin{split}
&\sum_{i=1}^p \int_M\left(A_{ii}^0(x_0) - 
\frac{n-2}{4(n-1)}S_g(x_0)\right)(\tilde u_\alpha^i)^2dv_g\\
&\le C \sum_{i\not= j}\int_{B_{x_0}(\delta)}\vert\tilde u_\alpha^i\vert \vert\tilde u_\alpha^j\vert dv_g +
C\left(\varepsilon + o(1)\right)\int_M\vert\tilde{\mathcal U}_\alpha\vert^2dv_g\\
&\hskip.4cm + \sum_{i=1}^p\left(\sup_{x \in B_{x_0}(\delta)}\left\vert A^0_{ii}(x) - A^0_{ii}(x_0)\right\vert\right)
\int_M\vert\tilde{\mathcal U}_\alpha\vert^2dv_g\\
&\le C\left(\varepsilon + \varepsilon_\delta + o(1)\right)\int_M\vert\tilde{\mathcal U}_\alpha\vert^2dv_g
\end{split}
\end{equation}
for all $\alpha$, where 
$o(1) \to 0$ as $\alpha \to +\infty$, $\varepsilon_\delta \to 0$ as $\delta \to 0$, 
and $C > 0$ does not depend on $\alpha$, $\varepsilon$, and $\delta$. By (\ref{IneqSecProof4Eqt1}) 
there exists $\varepsilon_0 > 0$ such that
\begin{equation}\label{Step6Sec4ProofEqt5}
A^0_{ii}(x_0) \ge \frac{n-2}{4(n-1)}S_g(x_0) + \varepsilon_0
\end{equation}
for all $i$. Then the contradiction easily follows from (\ref{Step6Sec4ProofEqt4}) by 
choosing $\varepsilon > 0$ and $\delta > 0$ sufficiently small such that 
$C(\varepsilon + \varepsilon_\delta) \le \varepsilon_0/2$, where $C > 0$ is the constant in 
(\ref{Step6Sec4ProofEqt4}), and $\varepsilon_0$ is as in (\ref{Step6Sec4ProofEqt5}). 
This proves that for $\tilde{\mathcal U}_\alpha$ and $\tilde{\mathcal U}^0$ as in 
(\ref{Sec4ProofEqt1}) and (\ref{Sec4ProofEqt4}), 
we necessarily have that $\tilde{\mathcal U}^0 \not\equiv 0$ when 
we assume that $n \ge 4$ and that (\ref{IneqSecProof4Eqt1}) holds. Then, by 
Step \ref{Step1ProofLemSec4}, we get that (\ref{MainResToProve}) is also true. By 
standard elliptic theory, as already mentionned, 
this implies in turn that Lemma \ref{TheLem} is true.
\end{proof}

A possible extension of Lemma \ref{TheLem} is to replace 
the condition in the Lemma that
$A^0_{ii}(x) > \frac{n-2}{4(n-1)}S_g(x)$
for all $i$ and all $x$, 
by the condition that for any $i$, either $A^0_{ii}(x) > \frac{n-2}{4(n-1)}S_g(x)$ for all $x$, or 
$A^0_{ii}(x) < \frac{n-2}{4(n-1)}S_g(x)$ for all $x$, and hence that for any $i$, and any $x$,
\begin{equation}\label{Sec4FinRemEqt1}
A^0_{ii}(x) \not= \frac{n-2}{4(n-1)}S_g(x)\hskip.1cm .
\end{equation}
If we assume that the convergence 
of the $A^\alpha_{ij}$'s to the $A^0_{ij}$'s is in $C^1(M)$, and that 
the manifold is conformally flat, we can prove, 
with the estimates we obtained in 
Steps \ref{Step1ProofLemSec4} to \ref{Step5ProofLemSec4}, this claim that 
Lemma \ref{TheLem} remains true if we only assume  
(\ref{Sec4FinRemEqt1}). The 
proof, based on the Pohozaev identity instead of (\ref{BOLoc}), is as follows.
We let $\tilde{\mathcal U}_\alpha$ and $\tilde{\mathcal U}^0$ be given by 
(\ref{Sec4ProofEqt1}) and (\ref{Sec4ProofEqt4}). 
We assume by contradiction that 
${\mathcal U}^0 \equiv 0$, and let 
$x_0$ be the limit of the centers of the $1$-bubble from which 
the $p$-bubble $({\mathcal B}_\alpha)_\alpha$ in (\ref{Step1StatemSec4Eqt1}) is defined. 
Since $g$ is conformally flat, 
there exist $\delta_0 > 0$ and a conformal metric $\hat g$ to $g$ such that $\hat g$ is flat in $B_{x_0}(4\delta_0)$. 
Let $\hat g = \varphi^{4/(n-2)}g$, where $\varphi$ is smooth and positive, and 
$\hat u_\alpha^i = \tilde u_\alpha^i\varphi^{-1}$ for all $\alpha$ and $i$. 
By conformal invariance of the conformal Laplacian, and by (\ref{Sec4ProofEqt2}), 
\begin{equation}\label{Sec4FinRemEqt2}
\Delta\hat u_\alpha^i + \sum_{j=1}^p\hat A^\alpha_{ij}\hat u_\alpha^j = (\hat u_\alpha^i)^{2^\star-1}
\end{equation}
in $B_{x_0}(4\delta_0)$ for all $i$ and all $\alpha$, where $\Delta = \Delta_{\hat g}$ is the Euclidean Laplacian, and
$\varphi^{2^\star-2}
\hat A^\alpha_{ij} = A^\alpha_{ij} - \frac{n-2}{4(n-1)}S_g\delta_{ij}$. 
The Pohozaev identity in 
the Euclidean space reads as
\begin{equation}\label{PohozaevIdent}
\begin{split}
&\int_\Omega(x^k\partial_ku)\Delta u dx + \frac{n-2}{2} \int_\Omega u(\Delta u)dx\\
&= - \int_{\partial\Omega}(x^k\partial_ku)\partial_\nu ud\sigma 
+ \frac{1}{2} \int_{\partial\Omega}(x,\nu)\vert\nabla u\vert^2d\sigma \\
&\hskip.4cm 
- \frac{n-2}{2}\int_{\partial\Omega}u\partial_\nu ud\sigma\hskip.1cm ,
\end{split}
\end{equation}
where $\nu$ is the outward unit normal to $\partial\Omega$,  
$d\sigma$ is the Euclidean volume element on $\partial\Omega$, and 
there is a sum over $k$ from $1$ to $n$.
For $\delta > 0$ small, we let  $\eta$ be a smooth cutoff function such that $\eta = 1$ in 
$B_{x_0}(\delta)$, $\eta = 0$ in $M\backslash B_{x_0}(2\delta)$, and $0 \le \eta \le 1$. We plugg the 
$\eta\hat u_\alpha^i$'s into the Pohozaev identity (\ref{PohozaevIdent}) and sum over $i$. 
In the process, we regard the $\hat u_\alpha^i$'s, $\varphi$, $\eta$, and the $\hat A^\alpha_{ij}$'s as 
defined in the Euclidean space. Thanks to 
(\ref{Sec4FinRemEqt2}), to the $C^1$-convergence of the $A^\alpha_{ij}$'s to the $A^0_{ij}$'s, 
and to $L^2$-concentration, coming back to the manifold, we get after lengthy 
(but simple) computations that
\begin{equation}\label{Sec4FinRemEqt3}
\begin{split}
&\sum_{i,j=1}^p \int_{B_{x_0}(\delta)}\left(A_{ij}^0(x_0) - 
\frac{n-2}{4(n-1)}S_g(x_0)\delta_{ij}\right)\tilde u_\alpha^i\tilde u_\alpha^jdv_g\\
&= \varepsilon_\delta O\left(\int_M\vert\tilde{\mathcal U}_\alpha\vert^2dv_g\right) 
+ o\left(\int_M\vert\tilde{\mathcal U}_\alpha\vert^2dv_g\right)
\end{split}
\end{equation}
for all $\alpha$, 
where $\varepsilon_\delta > 0$ is independent of $\alpha$ and such that 
$\varepsilon_\delta \to 0$ as $\delta \to 0$, and where the first term in the right hand side of 
(\ref{Sec4FinRemEqt3}) depends on $\delta$ only by $\varepsilon_\delta$. 
Let $i_0 = 1,\dots,p$, be such that 
${\mathcal B}_\alpha^{i_0} = B_\alpha$ for all $\alpha$, 
where $(B_\alpha)_\alpha$ is the $1$-bubble 
from which the $p$-bubble $({\mathcal B}_\alpha)_\alpha$ in (\ref{Step1StatemSec4Eqt1}) is defined. The argument 
we developed in the proof of Step \ref{Step5ProofLemSec4} gives that for $i \not= i_0$,
 \begin{equation}\label{Sec4FinRemEqt4}
\int_{B_{x_0}(\delta)}(\tilde u_\alpha^i)^2dv_g 
\le \varepsilon_\delta \int_M\vert\tilde{\mathcal U}_\alpha\vert^2dv_g
 \end{equation}
for all $\alpha$, where $\varepsilon_\delta > 0$ is independent of $\alpha$ and such that 
$\varepsilon_\delta \to 0$ as $\delta \to 0$. Combining the off diagonal estimates 
(\ref{Step5StatemSec4Eqt1}) of Step \ref{Step5ProofLemSec4}, $L^2$-concentration, (\ref{Sec4FinRemEqt3}), 
and (\ref{Sec4FinRemEqt4}), it follows that 
\begin{equation}\label{Sec4FinRemEqt5}
\begin{split}
&\left(A_{i_0i_0}^0(x_0) - 
\frac{n-2}{4(n-1)}S_g(x_0)\right)\int_M\vert\tilde{\mathcal U}_\alpha\vert^2dv_g\\
&= \varepsilon_\delta O\left(\int_M\vert\tilde{\mathcal U}_\alpha\vert^2dv_g\right) 
+ o\left(\int_M\vert\tilde{\mathcal U}_\alpha\vert^2dv_g\right)
\end{split}
\end{equation}
for all $\alpha$, where $\varepsilon_\delta > 0$ is independent of $\alpha$ and such that 
$\varepsilon_\delta \to 0$ as $\delta \to 0$, and where the first term in the right hand side of 
(\ref{Sec4FinRemEqt5}) depends on $\delta$ only by $\varepsilon_\delta$. 
Choosing $\delta > 0$ sufficiently small, and letting 
$\alpha \to +\infty$, we get a contradiction by combining (\ref{Sec4FinRemEqt1}) wth $i = i_0$ 
and (\ref{Sec4FinRemEqt5}). This proves that 
if we assume that the convergence of the $A^\alpha_{ij}$'s to the $A^0_{ij}$'s is in $C^1(M)$, and that 
the manifold is conformally flat, then 
Lemma \ref{TheLem} remains true if we replace the condition $A^0_{ii}(x) > \frac{n-2}{4(n-1)}S_g(x)$ 
for all $i$ and all $x$, by the condition that $A^0_{ii}(x) \not= \frac{n-2}{4(n-1)}S_g(x)$ 
for all $i$ and all $x$.

\medskip\noindent{\bf Acknowledgements:} The author is indebted to 
Olivier Druet and Fr\'ed\'eric Robert 
for their valuable comments on the manuscript.

\end{document}